\documentclass[11pt, bezier,amstex]{article}
\usepackage{amssymb,amsmath,latexsym}
\usepackage{epsfig}

\usepackage{epic,eepic,wrapfig,color}

\newtheorem{thm}{Theorem}[section]
\newtheorem{cor}[thm]{Corollary}
\newtheorem{lemma}[thm]{Lemma}
\newtheorem{prop}[thm]{Proposition}
\newtheorem{defn}[thm]{Definition}

\newtheorem{remark}[thm]{Remark}
\newtheorem{example}[thm]{Example}
\numberwithin{equation}{section}

\def\pf{{\medskip\noindent {\bf Proof. }}}
\def\qed{{\hfill $\Box$ \bigskip}}
\def\R{{\mathbb R}}
\def\P{{\mathbb P}}
\def\E{{\mathbb E}}
\def\1{{\bf 1}}

 \def\sB {{\cal B}}

 \def\bH {{\mathbb H}}

\def\bP {{\mathbb P}}  \def\bR {{\mathbb R}}

\def\R {{\mathbb R}}

\def\nn{\nonumber}

\def\wt{\widetilde}

\def\E{{\mathbb E}}
\def\P{{\mathbb P}}

\def\bea{\begin{align*}}
\def\eea{\end{align*}}
\def\bee{\begin{equation}}
\def\eee{\end{equation}}

\def\eps{\varepsilon}

\makeatletter
\@addtoreset{equation}{section}

\makeatother

\begin{document}
\bibliographystyle{plain}

\title{\Large \bf
Global Dirichlet Heat Kernel Estimates for
\\ Symmetric  L\'evy Processes  in Half-space}

\author{{\bf Zhen-Qing Chen}\thanks{Research partially supported
by   NSF Grant   DMS-1206276.}
\quad and \quad{\bf Panki Kim}\thanks{
This work was supported by the National Research Foundation of Korea(NRF) grant funded by the Korea government(MEST) (NRF- 2013R1A2A2A01004822)
}}
 %\date

\maketitle

\begin{abstract}
In this paper, we derive
explicit sharp two-sided estimates
 for the Dirichlet heat kernels of
a large class of symmetric
(but not necessarily rotationally symmetric)
L\'evy processes on  half spaces for all $t>0$.
These L\'evy processes may or may not have Gaussian component.
When L\'evy
density is comparable to a decreasing function with damping exponent $\beta$,
our estimate is explicit in terms of the distance to the boundary, the L\'evy exponent
and the damping exponent $\beta$ of L\'evy
density.
\end{abstract}

\bigskip
\noindent {\bf AMS 2000 Mathematics Subject Classification}: Primary
60J35, 47G20, 60J75; Secondary 47D07

\bigskip\noindent
{\bf Keywords and phrases}:
Dirichlet heat kernel, transition density,
survival probability, exit time, L\'evy system,
 L\'evy process, symmetric  L\'evy process.

\bigskip

\section{Introduction}

Classical Dirichlet heat kernel is  the fundamental solution of the heat equation
in
 an open set with zero boundary values.
Except for
a few special cases,
explicit form of the Dirichlet heat kernel is impossible to obtain.
Thus the best thing we can hope for is to establish  sharp two-sided estimates of Dirichlet heat kernels.
See \cite{D2} for upper bound estimates
and \cite{Zq3} for the lower bound estimate for Dirichlet heat kernels of diffusions in bounded $C^{1,1}$ domains.

The generator of a discontinuous L\'evy process is an integro-differential
operator and so it is a non-local operator.
Dirichlet heat kernels (if they exist) of
the generators of discontinuous L\'evy processes on an open set $D$
are the transition densities of such  L\'evy processes
killed upon leaving  $D$.
Due to this connection,
 obtaining sharp estimates on   Dirichlet heat kernels
  is a fundamental problem both
in probability theory and in analysis.

Before \cite{CKS}, sharp two-sided
estimates for the Dirichlet heat kernel of any non-local operator
in
open sets are unknown.
Jointly with R. Song,  in \cite{CKS} for the fractional Laplacian
$\Delta^{\alpha /2}:=-(-\Delta)^{\alpha /2}$ with zero
exterior condition, we succeeded in establishing sharp two-sided estimates in any $C^{1, 1}$
open set
$D$ and over any finite time interval
(see \cite{BGR} for an extension to non-smooth open sets).
When $D$ is bounded, one can easily deduce large time heat kernel estimates
from short time estimates by a spectral analysis.
The approach developed in \cite{CKS} provides a road map
for establishing sharp two-sided heat kernel estimates of
other discontinuous processes in open subsets of $\R^d$
(see \cite{BG, BGR, BGR3, CKS5, CKSa1, CKS7, CKS9, CKS6, KK}).
In \cite{CT, CKS4, CKS7},  sharp
two-sided estimates for the Dirichlet heat kernels $p_D(t, x, y)$ of
$\Delta^{\alpha/2}$ and of
$m-(m^{2/\alpha}-\Delta)^{\alpha/2}$ are obtained for all $t>0$ in two classes of unbounded open sets: half-space-like $C^{1, 1}$ open sets
and
exterior open sets. Since the estimates in \cite{CT, CKS4, CKS7} hold for all $t>0$, they are called global Dirichlet heat kernel estimates.
An important question in this direction is
for
how general discontinuous L\'evy processes one can prove sharp two-sided global Dirichlet heat kernel estimates
in unbounded open subsets of $\R^d$.

We conjectured in \cite[(1.9)]{CKS4}
that, when $D$ is a half space-like $C^{1,1}$ open set,
the following two-sided estimates
hold for a large class of rotationally symmetric
L\'evy process $X$ whose L\'evy exponent of $X$ is
$\Psi(|\xi|)$:
there are constants $c_1,
c_2, c_3\geq 1$ such that for every $(t, x, y)\in (0, \infty) \times D\times
D$,
\begin{eqnarray}
&& \hskip -0.6truein \frac1{c_1} \left(\frac1{\sqrt{t\, \Psi (1/\delta_D(x))}}
\wedge 1 \right)  \left(\frac1{\sqrt{t\, \Psi (1/\delta_D(y))}} \wedge 1 \right)
 p(t, c_2(y-x)) \,\le\, p_D(t, x,  y) \nonumber \\
&& \le\, c_1 \, \left(\frac1{\sqrt{t\, \Psi (1/\delta_D(x))}}
\wedge 1 \right)  \left(\frac1{\sqrt{t\, \Psi (1/\delta_D(y))}} \wedge 1\right)
 \, p(t, c_3(y-x)) \label{e:1.10}
\end{eqnarray}
where  $p(t, x)$ is the transition density of $X$.
In this paper, we use ``:=" as a way of definition.
For $a, b\in \bR$,  $a\wedge b:=\min \{a, b\}$
and $a\vee b:=\max\{a, b\}$.

Recently,  the above conjecture is confirmed in \cite[Theorem 5.8]{BGR3}  for   rotationally
symmetric unimodal L\'evy
process
whose L\'evy  exponent $x \to \Psi(|x|)$
satisfies the following upper and lower scaling properties: there are constants
$0<\beta_1<\beta_2 <2 $ and $C>1$ so that
\begin{equation}\label{e:1.2}
c^{-1} \left(\frac{R}{r} \right)^{\beta_1} \leq \frac{\Psi (R ) }{\Psi (r ) } \leq
c \left(\frac{R}{r} \right)^{\beta_2}
\quad \hbox{for any } R\geq r>0.
\end{equation}
Condition \eqref{e:1.2}
implies that the L\'evy process is  purely discontinuous, and by
\cite[Corollary 23]{BGR1},  its L\'evy intensity kernel
$x \to j(|x|)$ satisfies \begin{equation}
 \frac{c^{-1} }{|x|^d \Phi (|x|)} \le j(|x|)   \le  \frac{c }{|x|^d \Phi (|x|)}  \quad \text{for all }x\neq0,
\end{equation}
where $\Phi  ( r) :=  \max_{|x| \le r} 1/\Psi (1/|x|) $.
It is easy to see that $\Phi (r)$ is comparable to $1/\Psi (1/r)$.
It follows from \eqref{e:1.2}
that the
same two-sided estimates hold for $\Phi$ in place of $\Psi$.
Thus condition \eqref{e:1.2} excludes damped L\'evy  processes such as relativistic stable processes.
We remark here that under condition \eqref{e:1.2}, it follows
as a special case from \cite{CK2}
that the transition density $p(t, x)$ of the  rotationally symmetric unimodal L\'evy process has the following two-sided estimates:
\begin{equation}\label{e:1/4}
c^{-1} \left( \Phi^{-1}(t) \wedge \frac{t}{|x|^d \Phi (|x|)} \right) \leq p(t, x) \leq
c\left( \Phi^{-1}(t) \wedge \frac{t}{|x|^d \Phi (|x|)} \right)
\quad \hbox{for all } t>0 \hbox{ and } x\in \R^d.
\end{equation}

In this paper, we mainly focus
on estimate
 \eqref{e:1.10} when $D$ is a half space and we prove that \eqref{e:1.10}
holds for a large class of
symmetric L\'evy processes which may not be isotropic and may have damped L\'evy kernel.
Moreover, our
  symmetric L\'evy processes may or may not have Gaussian component.
Once the global Dirichlet heat kernel estimates  for upper half space
and short time heat kernel estimates  on $C^{1,1}$ open sets
 are obtained,
one can then use the ``push inward" method introduced in \cite{CT} to extend the results to half-space-like $C^{1,1}$ open sets.
See Remark \ref{R:7.2}.
For recent results on short time Dirichlet heat kernel estimates for symmetric L\'evy processes  in $C^{1,1}$ open sets,
we refer the reader to \cite{CKS9, CKS6}.
Note that, for all symmetric L\'evy process in $\R$, except compound  Poisson processes,
the survival probability $\P_x (\zeta >t)$
of  its subprocess in  the half line $(0, \infty)$
 is comparable to
$$ 1 \wedge \frac{ \max _{0\le y \le 1/x} 1/\sqrt{\Psi (y)}}{\sqrt{t}} ,
$$
  where $x \to \Psi(|x|)$ is its characteristic exponent
(see \cite[Theorem 4.6]{KMR} and   \cite[Proposition 2.6]{BGR2}).
This fact, which is used several times in this paper,
 is essential in our approach.

We now give more details on the main results of this paper.
In this paper, $d\geq 1$ and  $X=(X_t, \P_x)_{t \ge 0, x \in \R^d}$ is a symmetric discontinuous   L\'evy process  (but possibly with
Gaussian component) on $\R^d$  with L\'evy exponent
$\Psi(\xi)$ and L\'evy  density $J$
where $\P_x (X_0=x)=1$.
That is, $X$ is a right continuous symmetric
process having independent stationary increments with
\begin{align}\label{e:Psi}
\E_x\left[e^{i\xi\cdot(X_t-X_0)}\right]=e^{-t\Psi(\xi)}
\quad \quad \mbox{ for every } x\in \R^d \mbox{ and } \xi\in \R^d.
\end{align}
Throughout this paper, we assume that $X$ is not a compounded Poisson process,
which corresponds exactly to the case that $\Psi$ is unbounded.
It is known that
$$
\Psi(\xi)= \sum_{i, j=1}^d
a_{ij} \xi_i \xi_j+ \int_{\R^d}(1-\cos(\xi\cdot y))J(y)dy
\qquad \hbox{for } \xi=(\xi_1, \cdots, \xi_d)\in \R^d,
$$
where $A=(a_{ij})$ is a constant, symmetric, non-negative definite  matrix
and $J$ is a symmetric non-negative function on $\R^d\setminus \{0\}$ with
$\int_{\R^d} (1\wedge |z|^2) J(z) dz<\infty$.

When
$$
\int_{\R^d}\exp\left( -t \Psi(\xi)\right) d\xi < \infty \qquad \text{ for }t>0, \eqno({\bf ExpL})
$$
  the transition density
$p(t,x,y)=p(t,y-x)$
 of $X$ exists as a bounded
continuous function for each fixed $t>0$, and it is given by
$$
p(t, x):=(2\pi)^{-d}\int_{\R^d} e^{-i\xi \cdot x}e^{-t\Psi(\xi)}d\xi
, \quad t>0 .
$$
Moreover,
\begin{align}\label{e:asPip}
p(t,x) \leq(2\pi)^{-d}\int_{\R^d} e^{-t\Psi(\xi)}d\xi=p(t,0)<\infty.
\end{align}
Clearly,
condition ({\bf ExpL})
holds if $\inf_{|\xi| = 1} \sum_{i, j=1}^d a_{i j} \xi_i \xi_j>0$.
Conversely, suppose that, for every $t>0$,
 $X_t$ has a probability density function $p(t, x)$
under $\P_0$ and that $x\mapsto p(t, x)$ is $L^2$-integrable.  Then \eqref{e:Psi} can be rewritten as
$$ \int_{\R^d} p(t, x) e^{ix \cdot \xi} dx = e^{-t\Psi (\xi)}
$$
and so, by the Plancherel theorem,  $e^{-t\Psi (\xi)} \in L^2(\R^d)$ for every $t>0$; that is, ({\bf ExpL}) holds.
Hence condition ({\bf ExpL}) is equivalent to the existence of transition density function  $p(t, x)$ of $X$
that is $L^2(\R^d)$-integrable in $x$ for every $t>0$.
If we assume that $X$ has a transition density function $p(t, x)$ that is continuous in $x$,
then, since $\int_{\R^d} p(t, x)^2 dx = p(2t, 0)<\infty $,
  {\rm ({\bf ExpL})} holds.
See \cite[Proposition 4.1]{KS} for additional discussion
on condition ({\bf ExpL}).

Let
\bee \label{e:1.9Psi}\Psi^*(r):= \sup_{|z| \le r} \Psi(z)\eee
and use $\Phi$ to denote the non-decreasing  function
\bee \label{e:1.9}
\Phi(r)=\frac{1}{\Psi^*(1/r)}  \qquad \hbox{for } r>0  .
\eee
Note that since $\Psi (z)$ is a continuous unbounded function on $\R^d$,
$\Psi^*$ is a non-decreasing continuous function on $[0, \infty)$ with $\Psi^*(0)=0$ and
$\lim_{r\to \infty} \Psi^* (r)=\infty$. Consequently, 
 $\Phi$ is a non-decreasing continuous function on $[0, \infty)$ with $\Phi (0)=0$ and 
 $\lim_{r\to \infty} \Phi (r)=\infty$.
The right continuous  inverse function of $\Phi$ will
be denoted by the usual notation $\Phi^{-1}(r)$;
that is,
$$
    \Phi^{-1} (t)=\inf\{s>0: \Phi (s)>t\}.
$$
Note that $0<\Phi^{-1}(t)<\infty$ for every $t>0$ and $\lim_{t\to 0+} \Phi^{-1} (t)=0$.
Define for $r> 0$,
$$
\Psi_1^*(r):= \sup_{s \in (-r, r)} \Psi((0,\dots, 0, s)).
$$
We consider the following condition: there exists a constant $c \ge 1$ such that
$$
\Psi^*(r) \,\le \,c \,\Psi_1^*(r)  \quad \text{ for all } r>0.   \eqno({\bf Comp})
$$
Condition ({\bf Comp}) is a mild assumption that is satisfied by a large class of symmetric L\'evy processes, see Lemma \ref{l:Cnew}.
Under assumptions ({\bf ExpL}) and ({\bf Comp}), we derive in Lemma
\ref{l:u_off1p} a useful upper bound estimate for Dirichlet heat kernels.

In general, the explicit estimates of the  transition density $p(t,y)$ in $\R^d$
depend heavily on the
corresponding  L\'evy measure and Gaussian component (see \cite{CKK3, CK2}).
On the other hand,
 scale-invariant parabolic Harnack inequality holds with the explicit scaling in terms of L\'evy exponent for a large class of  symmetric L\'evy processes
(see  \cite[Theorem 4.12]{CK2}, \cite[Theorem 4.11]{CKK3} and our
Theorem \ref{T:1.1}).
Motivated by this, we first develop a rather general  version of
Dirichlet heat kernel upper bound estimate in Proposition \ref{p:u_off1n}
under the assumption
that  parabolic Harnack inequality  {\bf PHI($\Phi$)}  and ({\bf UJS}) hold.
See Section 3 for the definition of {\bf PHI($\Phi$)}.
We say ({\bf UJS}) holds if there exists a positive constant $c$ such that
for every $y\in \R^d$,
$$
 J(y) \le \frac{c}{r^d}
\int_{B(0,r)} J(y-z)dz \qquad
\hbox{whenever } r\le  |y|/2 . \eqno({\bf UJS})
 $$
 Note that
 ({\bf UJS}) is very mild assumption in our setting. In fact,
 ({\bf UJS}) always holds if $J(x) \asymp j(|x|)$ for
some non-increasing function $j$ (see  \cite[page 1070]{CKK2}). Moreover, if $J$ is continuous on $\R^d \setminus \{0\}$,
then ${\bf PHI(\Phi)}$  implies ({\bf UJS}).
In fact, using \eqref{e:ffgr} below instead of \cite[(2.10)]{CKK2},  this follows from the proof of \cite[Proposition 4.1]{CKK2}.
We also show in Theorem \ref{T:3.2} that ${\bf PHI(\Phi)}$  implies ({\bf ExpL}).

Assume in addition that for every $t>0$,  $x \to p(t, x)$ is weakly radially  decreasing in the following sense:
there exist constants $c>0$ and $C_1, C_2 >0$ such that
$$ p(t, x) \leq c  p(C_1t, C_2 y )
\qquad \hbox{for } t\in (0, \infty) \hbox{ and }
 |x|\geq |y|>0. \eqno({\bf HKC})
$$
We remark here that the same assumption with $C_1=1$ for small $t$ was made in \cite{CKS9}. Then our Dirichlet heat kernel upper bound estimate obtained in Propositions \ref{p:u_off1n}
yields the desired upper bound estimate in \eqref{e:1.10}.
Moreover,  we show that this assumption on $p(t, x)$
  (see Sections \ref{S:4} below) and the upper bound of $p_D(t, x)$ imply a very useful lower bound of $p_D(t, x)$;
see  Theorem  \ref{t:prlower}.

Jointly with T. Kumagai, in \cite{CK2, CKK3, CKK14}  we have established
 two-sided sharp heat kernel estimates for a large class of symmetric Markov
 processes.
 In Sections \ref{S:5}--\ref{S:7}, we assume the jumping kernels of our L\'evy process satisfy
the assumptions of \cite{CKK3, CKK14, CK2},
that is, conditions ({\bf UJS}), \eqref{e:5.4} and \eqref{e:5.5} of this paper.
Then all the aforementioned conditions ({\bf ExpL}), ({\bf Comp}),  {\bf PHI($\Phi$)},
 ({\bf HKC}) are satisfied.
Using the two-sided heat kernel estimates
for symmetric Markov processes on $\R^d$ from
\cite{CKK3, CKK14, CK2}
 (see Theorem \ref{T:1.1})   and our lower
bound estimates for Dirichlet heat kernels in Theorem \ref{t:prlower}, we obtain two-sided global Dirichlet
heat kernel estimates \eqref{e:7.4}, essentially
prove the conjecture \eqref{e:1.10} for such symmetric L\'evy
processes and for $D=\bH$. See Remark \ref{R:7.2}(i) for details.
Furthermore,
our estimates are explicit
in terms of the distance to the boundary, the L\'evy exponent
and the damping exponent $\beta$ of L\'evy
density;
see Theorem \ref{t:final}.

In this paper, we use the following notations.
For any two positive functions $f$ and $g$,
$f\asymp g$ means that there is a positive constant $c\geq 1$
so that $c^{-1}\, g \leq f \leq c\, g$ on their common domain of
definition.
For any open set $V$, we denote by $\delta_V (x)$
the distance of a point $x$ to the boundary of $V$, i.e.,
$\delta_V(x)=\text{dist} (x,\partial V)$. We sometimes write point $z=(z_1, \dots, z_d)\in \bR^d$ as
$(\wt z, z_d)$ with $\wt z \in \bR^{d-1}$. We denote ${\bH}:=\{x=(\wt x,x_d)\in\bR^d:x_d>0\}$
the upper half space.
For a set $W$ in $\R^d$, $\overline{W}$ and $|W|$ denotes the closure and the Lebesgue measure of $W$ in $\R^d$,  respectively.
 Throughout the rest of this paper, the positive constants
 $a_0, a_1, M_1, C_i$,
$i=0,1,2,\dots $, can be regarded as fixed.
In the statements of results and the proofs, the constants $c_i=c_i(a,b,c,\ldots)$, $i=0,1,2,  \dots$, denote generic constants depending on $a, b, c, \ldots$, whose exact values are unimportant.
They start
anew in each statement and each proof.
The dependence of the constants on the dimension $d \ge 1$
may not be mentioned explicitly.

\section{Setup and preliminary estimates}\label{S:int}

Let $X$ be a symmetric L\'evy process on $\R^d$ with L\'evy exponent $\Psi (z)$ and
L\'evy  density $J (z)$.
Recall the definition of the non-decreasing functions $\Psi^* ( r) $ and
 $\Phi ( r )$ from
\eqref{e:1.9Psi} and \eqref{e:1.9}, respectively.
We emphasize that the L\'evy process $X$ does not need to be rotationally symmetric.
The following is known and true for any
 negative definite function
 (see \cite[Lemma 1]{G}).

\begin{lemma}\label{L:1.1}
For every  $t>0$ and $\lambda \geq 1$,
$$
1 \le \frac{\Phi(\lambda t)}{\Phi(t)} \le 2(1+\lambda^2).
$$
\end{lemma}

 For an open set $D$,   denote by  $\tau_D :=\inf\{t>0: \, X_t\notin D\}$
 the first exit
time of $D$.

\begin{thm} \label{l:hku0}
There exists a constant $c=c(d)>0$ such that
\begin{equation} \label{e:2.1P}
\P_0 (|X_t|>r) \,\le\,  c\,  t/ \Phi( r)  \quad
\hbox{for } (t, r) \in (0, \infty)\times (0, \infty).
\end{equation}
Consequently,   there exists   $\eps_1=\eps_1 (d) >0$ such that for all  $r >0$,
\begin{equation} \label{e:2.1}
 \P_0 \left( \tau_{B(0, r/2 )} >\eps_1 \Phi (r) \right) \ge 1/2.
\end{equation}
\end{thm}

 \pf \eqref{e:2.1P} is a consequence of \cite[(3.2)]{P} and  \cite[Corollary 1]{G}.

Since the L\'evy process $X$ is conservative, \eqref{e:2.1P} implies by \cite[Lemma 3.8]{BBCK} that for every $t, r>0$,
$$
\P_0 \left( \tau_{B(0, 2r)}\leq t \right)
=\P_0 \left( \sup_{s\leq t} |X_s|>2r \right) \leq 2 c_2 t /\Phi (r).
$$
Thus
$\P_0 \left( \tau_{B(0,  r/2)}\leq \eps_1 \Phi (r) \right)
\leq  2c_2 \eps_1 \Phi (r)/\Phi (r/2)$,
which by Lemma \ref{L:1.1} is no larger than  $20 c_2\eps_1$.
Taking $\eps_1=1/(40c_2)$ proves the theorem.
\qed

Recall that $J$ is the L\'evy density of $X$, which  gives rise to a L\'evy system for $X$ describing the jumps of  $X$.
For any $x\in \bR^d$, stopping time $S$ (with respect to the filtration of $X$), and nonnegative measurable function $f$ on $\bR_+ \times \bR^d\times \bR^d$ with $f(s, y, y)=0$ for all $y\in\bR^d$ and $s\ge 0$ we have
\begin{equation}\label{e:levy}
\E_x \left[\sum_{s\le S} f(s,X_{s-}, X_s) \right] = \E_x \left[ \int_0^S \left(\int_{\bR^d} f(s,X_s, y) J(X_s-y) dy \right) ds \right]
\end{equation}
(e.g., see \cite[Appendix A]{CK2} and the proof of  \cite[Lemma 4.7]{CK1}).

The following is a special case of \cite[Corollary 1]{G},
whose upper bound will be used in the next lemma.

\begin{lemma}\label{l:j-upper} For every   $r > 0$,
$$
\frac{1}{2\Phi(r)}\,  \le \,  \frac{\|A\|}{ r^2}+  \int_{\R^d}
  J(z) \left( 1\wedge \frac{|z|^2}{r^2} \right) dz
    \, \le \, \frac{8(1+2d)}{\Phi(r)}
$$
 where
$$
\|A\|: =
\sup_{|\xi| \le 1} \sum_{i, j=1}^d a_{i, j} \xi_i \xi_j.
$$
 \end{lemma}

 Using Lemma \ref{l:j-upper}, the proof of the next lemma is rather routine (see \cite[Lemma 4.10]{KSV7}). In fact,
this lemma is proved in \cite[Lemma 3 and Corollary 1]{G} for $a=1/2$. The proof for general $a$ is similar.
Thus we skip the proof.

\begin{lemma}\label{l2.1}
For every $a \in (0, 1)$, there exists $c=c(a)>0$ so that
for any $r > 0 $ and any open set $U$ with $U\subset B(0, r)$,
$$
{\P}_x\left(X_{\tau_U} \in B(0, r)^c\right) \,\le\, \frac{c}{ \Phi(r)}\E_x[\tau_U], \qquad x \in U\cap B(0, ar)\, .
$$
\end{lemma}

Note that for
$d$-th coordinate $X^d_t$ of
$X_t=(X^1_t, \dots, X^d_t)$
is a L\'evy process with
$$
\E_x\left[e^{i\eta (X^d_t-X^d_0)}\right]=\E_{(\wt 0, x)} \left[e^{i(\wt 0,\eta)\cdot(X_t-X_0)}\right] =e^{-t\Psi((\wt 0,\eta))}
\qquad  \mbox{ for every } x\in \R  \mbox{ and } \eta\in \R.
$$
That is,
$X_t^d$ is  a $1$-dimensional
symmetric L\'evy process  with L\'evy exponent
$\Psi_1(\eta) :=\Psi((\wt 0,\eta))$.
Throughout this paper we let $\Psi_1^*(r):= \sup_{z \in (-r, r)} \Psi_1(z)$ and use $\Phi_1$ to denote the increasing function
$$
\Phi_1(r)=\frac{1}{\Psi_1^*(r^{-1})}, \qquad r>0.
$$
Clearly
$$
\Psi_1^*(r) \le \Psi^*(r)\quad \text{ and }\quad   \quad \Phi(r) \le \Phi_1(r).
$$

Since $\tau_{\bH}:=\inf\{t>0: X_t^d >0\}$, by
\cite[Proposition 2.6]{BGR2} (see also \cite[Theorems 3.1 and 4.6]{KMR})
 all symmetric L\'evy processes, except compound Poisson processes, enjoy the following estimates of the survival probability on $\bH$.

\begin{lemma}\label{l:spH}
 Suppose that $\Psi_1$ is unbounded, then there exists
an absolute constant $C>0$ such that
$$
C^{-1}
\left(\sqrt{\frac{\Phi_1(\delta_{\bH}(x))}{t}} \wedge 1\right) \le \P_{x}(\tau_{\bH}>t) \,\le\, C\
\left(\sqrt{\frac{\Phi_1(\delta_{\bH}(x))}{t}} \wedge 1\right) .
$$
\end{lemma}

Let
$$
\tau^1_{r}:= \inf\left\{t>0: X^d_t \notin \big(0, r \big)\right\}
$$
Combining  \cite[Lemma 2.3 and Proposition 2.4]{BGR2}  we have

\begin{lemma}\label{l:tau}
Suppose that $\Psi_1$ is unbounded, then there exists
 an absolute constant
 $c>0$
 such that for any  $r\in (0, \infty)$ and
\begin{eqnarray*}
\E_{(\wt 0, x)} [\tau^1_{r}]\,\le \, c
\,  \Phi_1(r)^{1/2} \Phi_1(\delta_{(0, r)}(x))^{1/2}
\quad  \hbox{ for }  x\in (0,  r).
\end{eqnarray*}
\end{lemma}

\medskip

Recall that, when ({\bf ExpL}) holds,
  the transition density
$p(t, x, y):=p(t,x-y)$
of $X$ exists as a bounded
continuous function.
In this case, for an open set $D$ we define
\begin{equation}\label{e:hkos1}
p_D(t,x,y)  :=  p(t,x,y) - \E_x [ \>p(t - \tau_D,
X_{\tau_D},y) : \tau_D < t]\quad \text{ for }
t>0, x,y \in D
\end{equation}
Using the strong Markov property of $X$, it is easy to verify that
$p_D(t,x,y)$ is the transition density for $X^D$,
the subprocess of $X$ killed upon leaving an open set $D$.

\begin{lemma}\label{l:u_near0}
Suppose {\rm ({\bf ExpL})} holds. Then for every  $(t,x,y) \in (0, \infty) \times \bH
\times \bH$,
$$
 p_{\bH}(t,x,y) \,\le\,
 3C^2p(t/3,0)
\left(\sqrt{\frac{\Phi_1(\delta_{\bH}(x))}{t}} \wedge 1\right)\left(\sqrt{\frac{\Phi_1(\delta_{\bH}(y))}{t}} \wedge 1\right)
$$
where
$C$ is the constant in Lemma \ref{l:spH}.
\end{lemma}

\pf Since by \eqref{e:asPip}
$$
\sup_{z,w \in  \bH} p_{\bH} (t/3,z,w)\le \sup_{z \in  \bH} p (t/3,z) =p(t/3,0),$$
using the semigroup property and symmetry we have
\begin{align*}
 p_{\bH}(t,x,y) &=\int_{\bH} \int_{\bH} p_{\bH}(t/3,x,z)
p_{\bH}(t/3,z,w) p_{\bH}(t/3,w,y) dzdw\\
& \le p(t/3,0) \, \P_{x}(\tau_{\bH}>t/3)
\P_{y}(\tau_{\bH}>t/3).
\end{align*}
Now the lemma follows from
Lemma \ref{l:spH}.
\qed

Using \eqref{e:levy}, the proof
of next lemma is the same as the one  in
\cite[Lemma 3.1]{CKS6} so
it is omitted.

\begin{lemma}\label{L:4.1}
Suppose {\rm ({\bf ExpL})} holds.
Suppose that $U_1, U_3, E$ are open subsets of $\mathbb{R}^d$, with $U_1, U_3\subset E$ and $dist (U_1, U_3)>0$.
Let $U_2:=E\backslash(U_1\cup U_3)$.
If $x\in U_1$ and $y\in U_3$, then for every $t>0$ we have
\begin{eqnarray}
p_E(t, x, y) &\le&  \mathbb{P}_x \left(X_{\tau_{U_1}}\in U_2\right)\cdot\sup _{s<t, z\in U_2} p_E(s, z, y)\nn\\
&& + \int_0^t \mathbb{P}_x\left(\tau_{U_1}>s\right) \mathbb{P}_y \left(\tau_E>t-s \right)ds \cdot\sup_{u\in U_1, z\in U_3} J(u- z) \label{eq:ub1}.
\end{eqnarray}
\end{lemma}

Recall condition ({\bf Comp}) from
the Introduction.  The next lemma says that it is a mild assumption.

\begin{lemma}\label{l:Cnew}
Suppose there are a non-negative  function $j$ on $(0, \infty)$ and $a\ge0$
$c_i \geq 1$, $i=1,2$, such that
\begin{equation}\label{e:psi2n}
c_1^{-1}a |y|^2 \leq \sum_{i, j=1}^d a_{i, j} y_i y_j \leq c_1a |y|^2 \quad\text{and}\quad
c_1^{-1} j(|y|/c_2)\le J(y) \le c_1 j(c_2|y|)  \quad \mbox{for all } y\in \R^d,
\end{equation}
Then $c^{-1} \Psi_1^* (r) \leq \Psi^* (r) \leq c\Psi_1^* (r)$, and so
 ${\rm ({\bf Comp})}$ holds.
\end{lemma}

\pf
Let
$$
\phi(|\xi|)=a|\xi|^2+ \int_{\R^d}(1-\cos(\xi\cdot y))j(|y|)dy.
$$
By a change of variables, \eqref{e:psi2n} implies that
\begin{align*}
\Psi(\xi)
& \le c_1\left(a|\xi|^2+ \int_{\R^d}(1-\cos(\xi\cdot y))j(c_2|y|)dy \right)\\
&\le c_3\left(a|c_2^{-1}\xi|^2+ \int_{\R^d}(1-\cos(c_2^{-1}\xi\cdot z))j(|z|)dz \right) = c_3\phi(|\xi|/c_2).
\end{align*}
and
$$
\Psi(\xi)
 \ge c_1^{-1}\left(a|\xi|^2+ \int_{\R^d}(1-\cos(\xi\cdot y))j(
 |y|/c_2)dy \right)
\ge  c_5\phi(c_2 |\xi|)
$$
Thus by Lemma \ref{L:1.1}, which holds for any  negative definite function,
 $\Psi^*(r)\asymp  \sup_{s \le r}\phi (s) \asymp \Psi_1^*(r)$ for all $r>0$.
 \qed

Using Lemma \ref{L:4.1},
 we can obtain the following upper bound of $p_{\bH}(t,x,y)$.

\begin{lemma}\label{l:u_off1p}
Suppose
{\rm ({\bf ExpL})} and {\rm ({\bf Comp})} hold.
 For each $a>0$, there exists a constant
 $c=c(a, \Psi)>0$ such that  for every $(t,x,y) \in (0, \infty)
\times \bH \times \bH$ with $  a\Phi^{-1}(t) \le |x-y|$,
\begin{align}
 p_{\bH}(t,x,y) \nn
 \le & c
\Big(\sqrt{\frac{\Phi(\delta_{\bH}(x))}{t}} \wedge 1\Big)
\Big(   \sup_{(s,z):s\le t,
 \frac{|x-y|}{2}  \leq |z-y| \le \frac{3 |x-y|}{2}}  p_{\bH}(s, z, y)
\\ &\qquad  \qquad      + \Big(\sqrt{t \Phi(\delta_{\bH}(y))} \wedge  t\Big)
\sup_{w:
 |w| \ge \frac{|x-y|}{3}}
J(w)\Big). \label{e:dfnewp}
\end{align}
\end{lemma}

\pf If $\delta_{\bH}(x) >  a\Phi^{-1}(t)/(24)$,
by Lemma \ref{L:1.1}
$$
\sqrt{\frac{\Phi(\delta_{\bH}(x))}{t}} \ge \sqrt{\frac{\Phi( a\Phi^{-1}(t)/(24))}{\Phi(\Phi^{-1}(t))}}
\ge \sqrt{\frac12 \frac{ a^2}{a^2+ (24)^2}}.
$$
Thus \eqref{e:dfnewp} is clear.

We now assume  $\delta_{\bH}(x) \le  a\Phi^{-1}(t)/(24) \le |x-y|/(24)$ and let $x_0 =(\wt x, 0)$, $U_1:=B( x_0,
 a\Phi^{-1}(t)/(12)) \cap \bH$, $U_3:= \{z\in \bH: |z-x|>|x-y|/2\}$
and $U_2:=\bH\setminus (U_1\cup U_3)$.
Recall that
$X_t^d$ is the $d$-th coordinate process of $X$  with L\'evy exponent
$\Psi_1(\eta) =\Psi((\wt 0,\eta))$.
Clearly,
 $$
 \tau_{U_1}\leq \inf\left\{t>0: X^d_t \notin \big(0,  a\Phi^{-1}(t)/12 \big)\right\}=:
 \tau^d_1.
 $$
Applying Lemma \ref{l:tau}
on the interval $\big(0,  a\Phi^{-1}(t)/12 \big)$ and assumption
{\rm ({\bf Comp})}, and noting
Lemma \ref{L:1.1}, we have
\begin{equation}\label{e:upae1p}
\E_{x}\left[ \tau_{U_1} \right]
\,\le\, \E^{X^d}_{\delta_\bH (x)} \big[ \tau^d_1\big] \leq
c_1\,  \sqrt{t \Phi(\delta_{\bH}(x))}.
\end{equation}
Since
$
|z-x| > 2^{-1}|x-y|  \ge a 2^{-1} \Phi^{-1}(t)$ for
$z\in U_3$,
we have for  $u\in U_1$ and $z\in U_3$,
$$
|u-z| \ge |z-x|-|x_0-x|-|x_0-u| \ge  \frac{1}{2}|x-y|- 6^{-1} a\Phi^{-1}(t) \ge \frac{1}{3}|x-y|.
$$
 Thus, $U_1 \cap U_3 = \emptyset$ and,
\begin{eqnarray}
\sup_{u\in U_1,\, z\in U_3}J(u-z) \le \sup_{(u,z):|u-z| \ge
\frac{1}{3}|x-y|}J(u-z)
= \sup_{w:|w|   \geq  \frac{1}{3}|x-y|}J(w).
\label{e:n01p}
\end{eqnarray}
Since for $z \in U_2$
$$
\frac32 |x-y| \ge |x-y| +|x-z| \ge  |z-y| \ge |x-y| -|x-z| \ge
\frac{|x-y|}2 \ge a2^{-1} \Phi^{-1}(t),
$$
we have
\begin{align}
&\sup_{s\le t,\, z\in U_2} p_{\bH}(s, z, y)
\le \sup_{ s\le t,  \frac{|x-y|}{2}  \leq |z-y| \le \frac{3 |x-y|}{2}}  p_{\bH}(s, z, y).\label{e:n02p}
\end{align}
Moreover, by Lemma \ref{l:spH} and {\rm ({\bf Comp})}
\begin{eqnarray*}
&& \int_0^t \mathbb{P}_x\left(\tau_{U_1}>s\right) \mathbb{P}_y \left(\tau_{\bH}>t-s \right)ds \le
\int_0^t \mathbb{P}_x\left(\tau_{\bH}>s\right) \mathbb{P}_y \left(\tau_{\bH}>t-s \right)ds \\
&\le & c_3 \ \int_0^t  \sqrt{\frac{\Phi_1 (\delta_{\bH}(y))}{s}} \left(\sqrt{\frac{\Phi_1 (\delta_{\bH}(y))}{t-s}} \wedge 1\right) ds\\
&\le & c_4 \sqrt{\Phi(\delta_{\bH}(x))} \left(\sqrt{\Phi(\delta_{\bH}(y))}  \wedge \sqrt{t}\right)
\int_0^t \frac1{\sqrt{s(t-s)}} ds \\
&= & c_5 \sqrt{\Phi(\delta_{\bH}(x))} \left(\sqrt{\Phi(\delta_{\bH}(y))}  \wedge \sqrt{t}\right).
\end{eqnarray*}
Applying this and \eqref{eq:ub1}, \eqref{e:upae1p}, \eqref{e:n01p}  and
\eqref{e:n02p}, we obtain,
\begin{eqnarray*}
p_{\bH}(t, x, y)
&\le& c_{6} \int_0^t \mathbb{P}_x\left(\tau_{U_1}>s\right) \mathbb{P}_y \left(\tau_{\bH}>t-s \right)ds
\sup_{w:|w| \ge
\frac{1}{3}|x-y|}J(w)\\
&&\quad +
c_{6} \P_x\Big (X_{\tau_{U_1}}\in U_2 \Big) \sup_{s\le t,\, z\in U_2} p(s, z, y) \\
&\le& c_{7} \sqrt{\Phi(\delta_{\bH}(x))} \left(\sqrt{\Phi(\delta_{\bH}(y))}  \wedge \sqrt{t}\right) \sup_{w:|w| \ge
\frac{1}{3}|x-y|}J(w)\\
&& \quad  +
c_{6}\P_x\Big(X_{\tau_{U_1}}\in U_2\Big) \sup_{ s\le t, \frac{|x-y|}{2}  \leq |z-y| \le \frac{3 |x-y|}{2}}  p(s, z, y).
\end{eqnarray*}
Finally,
applying Lemmas \ref{l2.1} and \ref{L:1.1}
 and then \eqref{e:upae1p}, we have
$$
\P_x \Big(X_{\tau_{U_1}}\in U_2 \Big) \le \P_x
\Big(X_{\tau_{U_1}}\in B( x_0, a\Phi^{-1}(t)/(12))^c\Big)
\,\le\,  \frac{c_8}{t}\, \E_{x}[\tau_{U_1}]
\le c_9  t^{-1/2} \,
\sqrt{\Phi(\delta_{\bH}(x))} .
$$
Thus we  have proved
\eqref{e:dfnewp}.
\qed

\begin{example}\label{ex:2}
Let $d=1$, and
$$
 j(y)=|y|^{-(1+\alpha)}\left(1 + \sum_{n=1}^\infty  n {\bf 1}_{[n, n +2^{-n}]}(|y|)\right),
$$
or
$$
j(y)=|y|^{-(1+\alpha)} + \sum_{n=1}^\infty  n {\bf 1}_{[n, n +2^{-n}]}(|y|).
$$
In either case, the L\'evy exponent $\Psi (\xi)$ for
symmetric L\'evy process having $j(y)$ as its L\'evy intensity
 is comparable to $|\xi|^\alpha$.
 So conditions {\rm ({\bf ExpL})} and {\rm ({\bf Comp})}  are satisfied. Consequently results in
 this section are valid for this L\'evy process.
However,  in either case,
 $j$  does not satisfy
{\rm ({\bf UJS})} at points $n$
when $n$ is large (nor \eqref{e:5.5} below).
\end{example}

\section{Consequences of parabolic Harnack inequality}

 Let $Z_s:=(V_s, X_s)$ be
the space-time process of $X$, where
$V_s=V_0- s$.
The law of the space-time process $s\mapsto Z_s$ starting from
$(t, x)$ will be denoted as $\mathbb{P}^{(t, x)}$.

\begin{defn}
A non-negative Borel  measurable function
$h(t,x)$ on $[0, \infty)\times \bR^d$ is said to be {\it parabolic}
(or {\it caloric})
on $(a,b]\times B(x_0,r)$
if for every relatively compact open subset $U$ of $(a,b]\times B(x_0,r)$,
$h(t, x)=
\mathbb{E}_{(t,x)}
 [h (Z_{\tau^Z_{U}})]$
for every $(t, x)\in U\cap ([0,\infty)\times
  \bR^d)$,
where
$\tau^Z_{U}:=\inf\{s> 0: \, Z_s\notin U\}$.
\end{defn}

It follows from the strong Markov property of $X$ and \eqref{e:hkos1}, $(t, x)\mapsto p_D(t,x,y)$ is parabolic on $(0, \infty) \times D$ for every $y\in D$.

Throughout this section, we assume
 the following (scale-invariant)
parabolic Harnack inequality
{\bf PHI($\Phi$)} holds
for $X$:
  For every $\delta\in (0, 1)$, there exists
$c=c(d, \delta)>0$
such that for every $x_0\in
\R^d$, $t_0\ge 0$, $R>  0$ and every non-negative function $u$ on $[0,
\infty)\times \R^d$ that is
  parabolic
 on $(t_0,t_0+4\delta
\Phi(R)]\times B(x_0,
R)$,
$$
\sup_{(t_1,y_1)\in Q_-}u(t_1,y_1)\le c \, \inf_{(t_2,y_2)\in
Q_+}u(t_2,y_2),
\eqno {\bf PHI(\Phi)}
$$
where $Q_-=(t_0+\delta \Phi(R),t_0+2\delta  \Phi(R)]\times B(x_0,
R/2)$
and $Q_+=[t_0+3\delta \Phi(R),t_0+ 4\delta \Phi(R)]\times B(x_0,
R/2)$.

\begin{thm}\label{T:3.2}
 Suppose that ${\bf PHI(\Phi)}$ holds. Then {\rm ({\bf ExpL})} holds  and so the L\'evy process $X$ has a
bounded continuous density function $p(t, x)$. Moreover,
  there is a constant $c >0$ so that
 \begin{equation}\label{e:ffgr}
 p(t, x)
\le     c   (\Phi^{-1}(t))^{-d} \qquad \hbox{for every } t>0 \hbox{ and } x\in \R^d.
\end{equation}
\end{thm}

\pf  Let $f\geq 0$ be an arbitary bounded $L^1$-integrable function on $\R^d$. Clearly $u(t, x): =P_t f (x)=\E_x [ f(X_t)]$
is a non-negative parabolic function on $(0, \infty)\times \R^d$.  Thus by ${\bf PHI(\Phi)}$ and the symmetry of the semigroup
$\{P_t; t>0\}$,  for every $x_0\in \R^d$
and $t>0$,
\begin{eqnarray*}
P_t f (x_0)
&\le&  c_1 \inf_{z\in B(x_0, \Phi^{-1}(t))} P_{3t} f (z)  \le c_2 (\Phi^{-1}(t))^{-d} \int_{B(x_0, \Phi^{-1}(t)) }   P_{3t} f(z) dz \\
&=& c_2 (\Phi^{-1}(t))^{-d} \int_{\R^d}  P_{3t} {\bf 1}_{    B(x_0, \Phi^{-1}(t))  }(z)   f(z) dz  \le  c_2 (\Phi^{-1}(t))^{-d} \int_{\R^d}   f(z) dz.
\end{eqnarray*}
This implies that $X$ has a transition density function $p(t, x)$ and
$p(t, x)\leq c_2  (\Phi^{-1}(t))^{-d}$
a.e. on $\R^d$.
Consequently, as mentioned earlier in the Introduction,  {\rm ({\bf ExpL})} holds by the Plancherel theorem,
which in turn implies that $p(t, x)$ is bounded and
  jointly continuous and so  \eqref{e:ffgr} holds.
\qed

Under the assumptions  ${\bf PHI(\Phi)}$ and  ({\bf UJS}),
we can
derive an interior lower bound for $p_D(t, x, y)$ for all $t>0$;
see Propositions \ref{step1} and \ref{step3}.
Similar bound for $t \le T$ was obtained in \cite{CKS6}
for subordinate Brownian motions with Gaussian component.
In this section, we use the convention that $\delta_{D}(\cdot) \equiv \infty$ when $D
=\bR^d$.

The next lemma holds for every symmetric L\'evy process and it
follows from \cite[(3.2)]{P} and \cite[Corollary 1]{G}.

\begin{lemma}\label{L:3.3}
 For any  positive
 constants  $a, b$, there exists
$c=c(a,b, \Psi)>0$ such that for all $z \in \bR^d$ and $t>0$,
$$
\inf_{y\in B(z, a \Phi^{-1}(t)/2)}\P_y \left(\tau_{B(z,
a \Phi^{-1}(t))}> bt  \right)\, \ge\, c.
$$
\end{lemma}

For the next two results,
$D$ is an arbitrary nonempty open set.

\begin{prop}\label{step1}
Suppose {\rm ${\bf PHI(\Phi)}$} holds.
Let $a>0$   be  a constant.
  There exists
 $c=c(a)>0$ such that
 \bee\label{e:lb1}
p_D(t,x,y) \,\ge\,c\, (\Phi^{-1}(t))^{-d}
\eee
for every $(t, x, y)\in (0, \infty)\times D\times D$ with
$\delta_D (x) \wedge \delta_D (y)
 \ge a \Phi^{-1}(t) \geq 4 |x-y|$.
 \end{prop}

\pf
We fix
$(t, x, y)\in (0, \infty)\times D\times D$ satisfying
$\delta_D(x)\wedge \delta_D(y)\geq a \Phi^{-1} (t) \geq 4 |x-y|$.
Note that
$|x-y| \le a\Phi^{-1}(t)/4 $
and that
$$
B(x, a\Phi^{-1}(t)/4) \subset B(y, a\Phi^{-1}(t)/2)\subset
 B(y, a\Phi^{-1}(t))  \subset D.
$$
So by the symmetry of $p_D$,
${\bf PHI(\Phi)}$, Theorem \ref{T:3.2},
and Lemma \ref{L:1.1},
 there exists $c_1=c_1(a)>0$ such that
$$
c_1 \, p_D(t/2, x, w) \, \le  \,  p_D(t,x,y) \quad \mbox{for every
} w \in B(x, a\Phi^{-1}(t)/4)  .
$$
This together with Lemma \ref{L:3.3}
yields that \begin{eqnarray*}
p_D(t, x, y) &\geq & \frac{c_1}{ | B(x, a\Phi^{-1}(t)/4)|}
\int_{B(x, a\Phi^{-1}(t)/4)} p_D(t/2, x, w)dw\\
&\geq & c_2 (\Phi^{-1}(t))^{-d} \, \int_{B(x, a \Phi^{-1}(t)/4)}
p_{B(x, a\Phi^{-1}(t)/4)} (t/2, x, w)dw \\
&=& c_2 (\Phi^{-1}(t))^{-d} \, \P_x \left( \tau_{B(x, a\Phi^{-1}(t)/4)} >
t/2\right) \,\geq \, c_3 \, (\Phi^{-1}(t))^{-d},
\end{eqnarray*}
where $c_i>0$ for $i=2, 3$.
\qed

Recall the condition  ({\bf UJS}) from the Introduction.

\begin{prop}\label{step3}
Suppose
${\bf PHI(\Phi)}$ and  {\rm({\bf UJS})} hold.
For every $a>0$,
there exists a constant $c=c(a)>0$ such that $p_D(t, x, y)\ge c t
 J(x-y)$ for every $(t, x, y)\in (0, \infty)\times D\times D$ with $\delta_D(x)\wedge \delta_D (y) \ge a \Phi^{-1}(t)$ and
$a \Phi^{-1}(t) \leq 4|x-y|$.
\end{prop}

\pf By Lemma~\ref{L:3.3}, starting at $z\in B(y, \,  (12)^{-1} a \Phi^{-1}(t))$, with probability at least $c_1=c_1(a)>0$ the process $X$ does not move more than $ (18)^{-1} a \Phi^{-1}(t) $ by time $t$.
Thus, using the strong Markov property and the L\'evy system in \eqref{e:levy}, we obtain
\begin{align}
&\P_x \left( X^D_t \in B \big( y, \,  6^{-1} a \Phi^{-1}(t) \big) \right)\nn\\
&\ge  c_1\P_x(X_{t\wedge \tau_{B(x,  (18)^{-1} a \Phi^{-1}(t))}}^D\in B(y,\,  (12)^{-1}a \Phi^{-1}(t))\hbox{ and }t \wedge \tau_{B(x,  (18)^{-1}a  \Phi^{-1}(t))} \hbox{ is a jumping time })\nonumber \\
&= c_1 \E_x \left[\int_0^{t\wedge \tau_{B(x,  (18)^{-1}a \Phi^{-1}(t))}} \int_{B(y, \,  (12)^{-1}a \Phi^{-1}(t))}
J(X_s- u) duds \right]. \label{e:nv1}
\end{align}
By ({\bf UJS}), we obtain
\begin{eqnarray}
 && \E_x \left[\int_0^{t\wedge \tau_{B(x,  (18)^{-1}a \Phi^{-1}(t))}} \int_{B(y, \,  (12)^{-1}a \Phi^{-1}(t))}
J(X_s-u) duds \right] \nn\\
 &= & \E_x \left[\int_0^{t} \int_{B(y, \,  (12)^{-1}a \Phi^{-1}(t))}
J(X^{B(x,  (18)^{-1}a \Phi^{-1}(t))}_s-u) duds \right] \nn\\
 &\ge & c_2\Phi^{-1}(t)^d \int_0^{t} \E_x \left[
J(X^{B(x,  (18)^{-1}a \Phi^{-1}(t))}_s-y) \right]ds \nn\\
& \ge & c_2\Phi^{-1}(t)^d \int_{t/2}^{t} \int_{B(x,  (72)^{-1}a \Phi^{-1}(t/2))}
J(w-y) p_{B(x,  (18)^{-1}a \Phi^{-1}(t))}(s, x,w)     dw ds.
   \label{e:nv2}
\end{eqnarray}
Since, for $t/2 <s <t$ and $w\in B(x,  (72)^{-1}a \Phi^{-1}(t/2))$
 $$ \delta_{B(x,  (18)^{-1}a \Phi^{-1}(t))} (w)
 \ge (18)^{-1}a \Phi^{-1}(t)-  (72)^{-1}a \Phi^{-1}(t/2)\ge 2^{-1}(18)^{-1}a \Phi^{-1}(s)
 $$
and
 $$|x-y| <  (72)^{-1}a \Phi^{-1}(t/2) \le 4^{-1}  (18)^{-1}a \Phi^{-1}(s),   $$
we have
by Theorem \ref{T:3.2} and
Lemma \ref{step1} that for $t/2 <s <t$ and $w\in B(x,  (72)^{-1}a \Phi^{-1} (t/2))$,
\begin{eqnarray}
p_{B(x,  (18)^{-1}a \Phi^{-1}(t))}(s, x,w)  \ge c_3\, (\Phi^{-1}(s))^{-d}  \ge c_3\, (\Phi^{-1}(t))^{-d}.
   \label{e:nv4}
\end{eqnarray}

Combining \eqref{e:nv1}, \eqref{e:nv2} with \eqref{e:nv4} and applying ({\bf UJS}) again, we get
\begin{align}
\P_x \left( X^D_t \in B \big( y, \,  6^{-1} a \Phi^{-1}(t) \big) \right)
 \ge & c_4 t  \int_{B(x,  (72)^{-1}a \Phi^{-1}(t/2))}
J(w-y)    dw
\nn\\
 \ge & c_5 t  (\Phi^{-1}(t/2))^{d}
J(x-y)   \ge  c_6 t  (\Phi^{-1}(t))^{d}
J(x-y).
   \label{e:nv6}
\end{align}
In the last inequality we have used Lemma \ref{L:1.1}.
The proposition now follows
 from the Chapman-Kolmogorov equation along with \eqref{e:nv1}, \eqref{e:nv2} and Proposition~\ref{step1}. Indeed,
\begin{eqnarray*}
p_D(t, x, y) &=& \int_{D} p_{D}(t/2, x, z)
p_{D}(t/2, z, y)dz\\
&\ge& \int_{B(y, \, a \Phi^{-1}(t/2)/ 6)}
p_{D}(t/2, x, z) p_{D}(t/2, z, y) dz\\
&\ge& c_7 (\Phi^{-1}(t/2))^{-d} \,
\bP_x \left( X^{D}_{t/2} \in B(y, a
\Phi^{-1}(t/2) /6) \right)\\
&\ge & c_6c_7\,{ t} J(x-y).
\end{eqnarray*}\qed

\medskip

We now apply Lemma \ref{l:u_off1p} to get the following heat kernel upper bound.

\begin{prop}\label{p:u_off1n}
Suppose   {\rm ({\bf Comp})}, ${\bf PHI(\Phi)}$ and  {\rm({\bf UJS})} hold.
Then there exists a constant
 $c>0$ such that  for every $(t,x,y) \in (0, \infty)
\times \bH \times \bH$
$$
 p_{\bH}(t,x,y) \le  c
\Big(\sqrt{\frac{\Phi(\delta_{\bH}(x))}{t}} \wedge 1\Big)\Big(\sqrt{\frac{\Phi(\delta_{\bH}(y))}{t}} \wedge 1\Big)
 \sup_{  |w| \ge  |x-y|/6 }  p(t, w) .
$$
\end{prop}

\pf By Lemma \ref{l:u_near0} and
Theorem \ref{T:3.2},
$$
 p_{\bH}(t,x,y) \,\le\, c_1 (\Phi^{-1}(t))^{-d}
\left(\sqrt{\frac{\Phi(\delta_{\bH}(x))}{t}} \wedge 1\right)\left(\sqrt{\frac{\Phi(\delta_{\bH}(y))}{t}} \wedge 1\right)  .
$$
If $\Phi^{-1}(t) \ge |x-y|$, by Proposition \ref{step1},
$ p(t,x-y) \ge  c_2  (\Phi^{-1}(t))^{-d}.$
Thus
\begin{align}
 p_{\bH}(t,x,y) \,\le\,c_3 \left(\sqrt{\frac{\Phi(\delta_{\bH}(x))}{t}} \wedge 1\right)\left(\sqrt{\frac{\Phi(\delta_{\bH}(y))}{t}} \wedge 1\right)  p(t,x-y). \label{e:dfghj1}
\end{align}

We extend the definition of
$p(t, w)$  by setting
 $p(t, w)=0$ for $t<0$
 and $w \in \R^d$.
 For each fixed $x, y\in \R^d$ and $t>0$ with $|x-y|> 8r$,
 one can easily check that $(s, w)\mapsto  p(s, w- y)$
is a parabolic function in $(-\infty, \infty)\times B(x, 2r)$.
Suppose  $  \Phi^{-1}(t) \le |x-y|$ and
let $(s,z)$ with $s\le t$ and  $\frac{|x-y|}{2}  \leq
 |z-y| \le \frac{3 |x-y|}{2}$.
 By ${\bf PHI(\Phi)}$,
 there is a constant
$c_4 \geq 1$
 so that for every $t>0$,
$$
\sup_{s \le t} p(s, z-y )\leq  c_4 p(t , z-y).
$$
Hence we have
\begin{equation}\label{eq:ppun31}
 \sup_{ s\le t,
 \frac{|x-y|}{2}  \leq |z-y| \le \frac{3 |x-y|}{2}}
p(s, z-y)  \leq c_4\sup_{
 \frac{|x-y|}{2}  \leq |z-y| \le \frac{3 |x-y|}{2}}
 p(t, z-y)=c_4\sup_{ \frac{|x-y|}{2}  \leq |z| \le \frac{3 |x-y|}{2}}  p(t, z).
 \end{equation}
Using this and Lemma \ref{l:u_off1p} and Proposition \ref{step3}, we have
 for every $(t,x,y) \in (0, \infty)
\times \bH \times \bH$ with $  \Phi^{-1}(t) \le |x-y|$,
\begin{align*}
& p_{\bH}(t,x,y) \\
\le & c_5
\Big(\sqrt{\frac{\Phi(\delta_{\bH}(x))}{t}} \wedge 1\Big)
\left(    \sup_{  \frac{|x-y|}{2}  \leq |z| \le \frac{3 |x-y|}{2}
 }  p(t, z)   +
 \Big(\sqrt{t \Phi(\delta_{\bH}(y))} \wedge  t\Big)
\sup_{|w| \ge \frac{|x-y|}{3}}
J(w)\Big)\right)\\
\le & c_6
\Big(\sqrt{\frac{\Phi(\delta_{\bH}(x))}{t}} \wedge 1\Big)
\left(  \sup_{|z| \ge |x-y|/2}  p(t, z)  + \sup_{
|w| \ge |x-y|/3 } p(t,w)\right)\\
\le & 2c_6
\Big(\sqrt{\frac{\Phi(\delta_{\bH}(x))}{t}} \wedge 1\Big)
 \sup_{|w| \ge  |x-y|/3}  p(t, w).
\end{align*}
In view of \eqref{e:dfghj1},
the last inequality holds in fact holds for all $(t,x,y) \in (0, \infty)
\times \bH \times \bH$.  Thus we have by an analogy of \eqref{eq:ppun31} that
for every $(t,x,y) \in (0, \infty)
\times \bH \times \bH$ with $|x-y|\geq \Phi^{-1}(t)$,
\begin{align}
 \sup_{  s \le t,  \frac{|x-y|}{2}  \leq |z-y| \le \frac{3 |x-y|}{2}}  p_{\bH}(s, z, y)
\le & c_7
\Big(\sqrt{\frac{\Phi(\delta_{\bH}(y))}{t}} \wedge 1\Big)
 \sup_{ |z-y|\ge |x-y|/2} \,
 \sup_{  |w|  \ge |z-y|/3}  p(t, w)\nn\\
 \le & c_{8}
\Big(\sqrt{\frac{\Phi(\delta_{\bH}(y))}{t}} \wedge 1\Big)
 \sup_{ |w| \ge  |x-y|/6 }  p(t, w). \label{e:dfghj11}
 \end{align}
Therefore by Lemma \ref{l:u_off1p}, Proposition \ref{step3} and
\eqref{e:dfghj11},
\begin{align*}
& p_{\bH}(t,x,y) \\
\le&
 c_{9}
\Big(\sqrt{\frac{\Phi(\delta_{\bH}(x))}{t}} \wedge 1\Big)
\left( \sup_{ s\le t,
 \frac{|x-y|}{2}  \leq |z-y| \le \frac{3 |x-y|}{2}}  p_{\bH}(s, z, y)   +
 \Big(\sqrt{t \Phi(\delta_{\bH}(y))} \wedge  t\Big)
\sup_{ |w| \ge \frac{|x-y|}{3}}
J(w)\Big)\right)\\
 \le&
 c_{10}
\Big(\sqrt{\frac{\Phi(\delta_{\bH}(x))}{t}} \wedge 1\Big)\Big(\sqrt{\frac{\Phi(\delta_{\bH}(y))}{t}} \wedge 1\Big)
\left( \sup_{|w| \ge |x-y|/6}  p(t, w)
+ \sup_{ |w| \ge |x-y|/3 } p(t,w)\right)\\
 \le&  2c_{10}
\Big(\sqrt{\frac{\Phi(\delta_{\bH}(x))}{t}} \wedge 1\Big)\Big(\sqrt{\frac{\Phi(\delta_{\bH}(y))}{t}} \wedge 1\Big)
 \sup_{
   |w| \ge |x-y|/6 }  p(t, w).
 \end{align*}
\qed

\section{Condition ({\bf HKC}) and its consequence}\label{S:4}

Under the condition ({\bf HKC}), clearly we have the following by Proposition \ref{p:u_off1n}.
\begin{thm}\label{p:u_off2n}
Suppose that conditions
{\rm ({\bf Comp})}, ${\bf PHI(\Phi)}$,  {\rm ({\bf HKC})},  and  {\rm({\bf UJS})} hold.
Then there exists a constant
 $C_3>0$ such that  for every $(t,x,y) \in (0, \infty)
\times \bH \times \bH$
$$
 p_{\bH}(t,x,y) \le C_3
\Big(\sqrt{\frac{\Phi(\delta_{\bH}(x))}{t}} \wedge 1\Big)\Big(\sqrt{\frac{\Phi(\delta_{\bH}(y))}{t}} \wedge 1\Big)
 p(C_1t,  6^{-1} {C_2}(x-y)).
$$
\end{thm}

We consider the following condition.
\begin{equation}\label{e:4.1}
\lim_{M \to  \infty} \sup_{r>0} \frac{\Psi^*( r)}{ \Psi^*( Mr )}  = 0.
 \end{equation}
It is equivalent to
\begin{align}\label{e:wlgPhi}
\lim_{M \to  \infty} \sup_{t>0} \frac{t}{ \Phi( M\Phi^{-1}(t))}   =\lim_{M \to  \infty} \sup_{t>0} t \Psi^*( M^{-1} {\Psi^*}^{-1}(1/t))= 0.
\end{align}

The following gives a sufficient condition for \eqref{e:4.1}.

\begin{prop}\label{p:ne12}
Suppose that $X$ has a transition density function
$p(t, x)$ that is continuous at  $x=0$ for every $t>0$ and $p(t,0) \le c  (\Phi^{-1}(t))^{-d}<\infty$ for all $t>0$.
Then condition
\eqref{e:4.1} holds. In particular,
${\bf PHI(\Phi)}$
implies that the condition
\eqref{e:4.1} holds.
\end{prop}

\pf
Since $X$ has a transition density function $p(t, x)$ that is continuous at $x=0$ for every $t>0$,
$\int_{\R^d} p(t/2, x)^2 dx = p(t, 0)<\infty$. It follows then $e^{-t\Psi (\xi)}$ is integrable on $\R^d$
and so
$$ p(t, x)= (2\pi)^{-d} \int_{\R^d} e^{-ix\cdot \xi} e^{-t \Psi (\xi )} d\xi.
$$
In particular, $ \int_{\R^d} e^{-t\Psi(\xi)}d\xi = (2\pi)^d p(t, 0) \leq  c_1  (\Phi^{-1}(t))^{-d}$
for every $t>0$. In other words,
$$
\int_{\R^d} e^{-\Psi(\xi)/r}d\xi  \le c_1  (\Phi^{-1}(1/r))^{-d}=c_1 ( (\Psi^*)^{-1}(r))^{d}, \quad r>0.
$$
Thus, for all $R, r>0$ we have
\begin{align}
\label{e:pnew1}
e^{-\Psi^*(R)/r} |B(0, R)| \le \int_{B(0, R)} e^{-\Psi(\xi)/r}d\xi  \le c_1 ( (\Psi^*)^{-1}(r))^{d} .
\end{align}
Note that $\Psi^*(r)$ is a non-decreasing continuous function on $[0, \infty)$ with
$\Psi^*(0)=0$ and \break
$\lim_{r\to \infty} \Psi^*(r)=\infty$.
Thus for every $r>0$ and $\lambda >1$, there is $R>0$ so that $\Psi^*(R)=\lambda r$.
Hence we have from \eqref{e:pnew1}  that
$
e^{-\lambda}   ( (\Psi^*)^{-1}(\lambda r))^{d} \le c_2 ( (\Psi^*)^{-1}(r))^{d} $, and so
\begin{align}
\label{e:pnew2}
\frac{ (\Psi^*)^{-1}(\lambda r)}{(\Psi^*)^{-1}(r)} \le (c_2 e^{\lambda})^{1/d} .
\end{align}
For $M>1$, let $\lambda=\lambda(M)=\log(M^d/c_2)$ so that $(c_2 e^{\lambda})^{1/d}=M$.
Then by \eqref{e:pnew2} with $s=(\Psi^*)^{-1}(r)$ we have
$Ms\ge(\Psi^*)^{-1}(\lambda r)$.
In other words,
$ \Psi^*(Ms) \ge \lambda r=\log(M^d/c_2)  \Psi^*(s)$.
Therefore
$$\sup_{s>0}\frac{\Psi^*( s)}{ \Psi^*( Ms )} \le  \frac{1}{\log(M^d/c_2)},$$ which goes to zero as
$M\to \infty$.
The last assertion of the theorem follows directly from Theorem \ref{T:3.2}.
\qed

 \begin{lemma}\label{l:ne12}
Suppose that
\eqref{e:4.1}
holds.
Then for each fixed $c>0$ the function
  $$ H_c(M):=
c^{-d} \sup_{t>0} \P \left(|X_t|>c M\Phi^{-1}(t) \right)
 $$
vanishes  at $\infty$;  that is, $\lim_{M \to  \infty} H_c(M)=0$.
\end{lemma}

\pf  By
Theorem \ref{l:hku0}, we have
$$
\sup_{t>0} \P \left(|X_t|> c M\Phi^{-1}(t) \right)
 \le  c_1\,   \sup_{t>0} \frac{t}{ \Phi( cM\Phi^{-1}(t))},
$$
which goes to zero as $M\to \infty$ by \eqref{e:wlgPhi}.
 \qed

For the remainder of this section, we
assume that  conditions
  {\rm ({\bf Comp})}, ${\bf PHI(\Phi)}$, {\rm ({\bf HKC})}   and  {\rm({\bf UJS})} hold,  and
discuss some lower bound estimates
of $p_{\bH}(t,x, y)$ under these conditions.
We first note that by {\rm ({\bf Comp})} and Lemma \ref{l:spH},
 there exists $C_0>0$ such that
\begin{align}\label{e:spH}
C_0^{-1}
\left(\sqrt{\frac{\Phi(\delta_{\bH}(x))}{t}} \wedge 1\right) \le \P_{x}(\tau_{\bH}>t) \,\le\, C_0\
\left(\sqrt{\frac{\Phi(\delta_{\bH}(x))}{t}} \wedge 1\right) .
\end{align}

We denotes by $e_d$
the unit vector in the positive direction of the $x_d$-axis in $\R^d$.
\begin{lemma}\label{l:loww_01}
There exist $a_1>0$ and $M_1 >4a_1$ such that for every $x \in \bH$  and $t>0$ we have
 $$\int_{ \{u \in {\bH} \cap B(\xi_x(t), M_1\Phi^{-1}(t)):\Phi(\delta_{\bH}(u)) >  a_1t\} }
 p_{{\bH}}(t,x,u)du
   \ge 4^{-1}  C_0^{-1} \left(\sqrt{\frac{\Phi(\delta_{\bH}(x))}{t}} \wedge 1\right)
   $$
   where $\xi_x(t):=x+a_1\Phi^{-1}(t)e_d$ and
   $C_0$ is the constant in
   \eqref{e:spH}.
\end{lemma}

\pf
By Theorem  \ref{p:u_off2n} and a change of variable, for every $t>0$ and $x \in \bH$,
\begin{align}\label{e:lbsp1}
&  \int_{
 \{ u \in \bH : \Phi(\delta_{\bH}(u)) \le  at\}
} p_{{\bH}}(t,x,u)du \nn\\
  \le & C_3\left(\sqrt{\frac{\Phi(\delta_{\bH}(x))}{t}} \wedge 1\right)
  \int_{
  \{ u \in \bH : \Phi(\delta_{\bH}(u)) \le  at\}} \left(\sqrt{\frac{\Phi(\delta_{\bH}(u))}{t}} \wedge 1\right) p(C_1 t,  6^{-1} {C_2}(x-u))
   du  \nn\\
     \le & C_3\sqrt{a}   \left(\sqrt{\frac{\Phi(\delta_{\bH}(x))}{t}} \wedge 1\right)
  \int_{
  \{ u \in \bH : \Phi(\delta_{\bH}(u)) \le  at\}} p(C_1 t,  6^{-1} {C_2}(x-u))  du  \nn\\
  \le &  C_3\sqrt{a}   \left(\sqrt{\frac{\Phi(\delta_{\bH}(x))}{t}} \wedge 1\right)  \int_{\R^d} p(C_1 t,  6^{-1} {C_2}(x-u))  du \nn\\
   =&  C_3(6/C_2)^{d}\sqrt{a}   \left(\sqrt{\frac{\Phi(\delta_{\bH}(x))}{t}} \wedge 1\right)  \int_{\R^d} p(C_1 t,  w)  dw
  =C_3(6/C_2)^{d}\sqrt{a}   \left(\sqrt{\frac{\Phi(\delta_{\bH}(x))}{t}} \wedge 1\right).
 \end{align}
Choose $a_1>0$ small so that
$C_3(6/C_2)^{d}\sqrt{a_1} \le (8 C_0)^{-1}$
where $C_0$ is the constant in
   \eqref{e:spH}.

For $x \in \bH$, we
let $\xi_x(t):=x+a_1\Phi^{-1}(t)e_d$.
For every $t>0$, ${M} \ge 2a_1 $ and $ u \in {\bH} \cap B(\xi_x(t), M\Phi^{-1}(t))^c$,  we have
$$|x-u| \ge |\xi_x(t)-u|-|x-\xi_x(t)| \ge |\xi_x(t)-u|  -a_1 \Phi^{-1}(t) \ge (1-\frac{a_1}{M}) |\xi_x(t)-u|\ge \frac12 |\xi_x(t)-u|$$
Thus using Theorem \ref{p:u_off2n} and  condition ({\bf HKC}), by a change of variable
 we have for every $t>0$ and  ${M} \ge 2a_1 $,
\begin{align}
& \int_{{\bH} \cap B(\xi_x(t), M\Phi^{-1}(t))^c} p_{{\bH}}(t,x,u)du
\nn \\
\le& C_3\left(\sqrt{\frac{\Phi(\delta_{\bH}(x))}{t}} \wedge 1\right)  \int_{{\bH} \cap B(\xi_x(t), M\Phi^{-1}(t))^c}\left(\sqrt{\frac{\Phi(\delta_{\bH}(u))}{t}} \wedge 1\right)
p(C_1 t,  6^{-1} {C_2}(x-u)) du  \nn\\
 \le& C_3\left(\sqrt{\frac{\Phi(\delta_{\bH}(x))}{t}} \wedge 1\right) \int_{{\bH} \cap B(\xi_x(t), M\Phi^{-1}(t))^c}p(C_1^2 t,  (12)^{-1} {C_2^2}(\xi_x(t)-u))  du  \nn\\
   \le& C_3\left(\sqrt{\frac{\Phi(\delta_{\bH}(x))}{t}} \wedge 1\right) \int_{B(0, M\Phi^{-1}(t))^c}p(C_1^2 t,  (12)^{-1} {C_2^2}u)
  du   \nn\\
    =& C_3((12)^{-1} {C_2^2})^d \left( \int_{B(0, (12)^{-1} {C_2^2} M\Phi^{-1}(t))^c}p(C_1^2 t,  v)
  dv \right)\left(\sqrt{\frac{\Phi(\delta_{\bH}(x))}{t}} \wedge 1\right)
  \nn\\
  \le &
 C_3 H_{(12)^{-1} {C_2^2}}(M)\left(\sqrt{\frac{\Phi(\delta_{\bH}(x))}{t}} \wedge 1\right).\label{n:neq11} \end{align}
 By Lemma \ref{l:ne12}, and Proposition \ref{p:ne12},
 we can
choose $M_1>4a_1$ large so that
$
C_3 H_{(12)^{-1} {C_2^2}}(M_1) < 8^{-1} \cdot C_0^{-1}.
$
Then by
   \eqref{e:spH}, \eqref{e:lbsp1}, \eqref{n:neq11} and our choice of $a_1$ and $M_1$, we conclude that
\begin{align}
 &\int_{
 \{u \in {\bH} \cap B(\xi_x(t), M_1\Phi^{-1}(t)):\Phi(\delta_{\bH}(u)) >  a_1t\}} p_{{\bH}}(t,x,u)du \nn\\
 =&
  \int_{{\bH}} p_{{\bH}}(t,x,u)du-\int_{{\bH} \cap B(\xi_x(t), M_1\Phi^{-1}(t))^c} p_{{\bH}}(t,x,u)du -\int_{ \{u \in {\bH}: \Phi(\delta_{\bH}(u)) \le   a_1t\}} p_{{\bH}}(t,x,u)du\nn \\
  \ge & 4^{-1} \cdot C_0^{-1} \left(\sqrt{\frac{\Phi(\delta_{\bH}(x))}{t}} \wedge 1\right).\nn
 \end{align}
\qed

For $x \in \bH$ and $t>0$,  let $\xi_x(t):=x+a_1\Phi^{-1}(t)e_d$ and define
\begin{align}\label{e:sB(x,t)}
\sB (x, t):= \left\{ z \in {\bH} \cap B(\xi_x(t), M_1\Phi^{-1}(t)):\Phi(\delta_{\bH}(z)) >  a_1t \right\}.
\end{align}

\begin{thm}\label{t:prlower}
There exist
$c_1, c_2>0$ such that for all $(t, x, y)\in (0, \infty) \times
\bH\times \bH$,
\begin{align}
p_{\bH}(t,x, y)
&\ge
 c_1  \left(\sqrt{\frac{\Phi(\delta_{\bH}(x))}{t}} \wedge 1\right)
\left(\sqrt{\frac{\Phi(\delta_{\bH}(y))}{t}} \wedge 1\right)\left(\inf_{(u,v) \in
\sB(x,t) \times \sB(y,t)}p_{\bH}(t/3,u,v) \right)
\label{e:loww1}  \\
&\ge
 c_2 \left(\sqrt{\frac{\Phi(\delta_{\bH}(x))}{t}} \wedge 1\right)
\left(\sqrt{\frac{\Phi(\delta_{\bH}(y))}{t}} \wedge 1\right)\times  \nn\\
& \qquad  \qquad \times
\begin{cases}
\displaystyle
\inf_{(u,v) : 2M_1 \Phi^{-1}(t) \le  |u-v|  \le 3|x-y|/2 \atop
\Phi(\delta_{\bH}(u)) \wedge \Phi(\delta_{\bH}(v)) >  a_1t}p_{\bH}(t/3,u,v) &
\hbox{ if }|x-y|  > 4M_1 \Phi^{-1}(t),\\
 (\Phi^{-1}(t))^d
& \hbox{ if }|x-y|   \le  4M_1 \Phi^{-1}(t).
\end{cases}
\label{e:loww21}
\end{align}
\end{thm}

\pf By Chapman-Kolmogorov equation,
\begin{align*}
&p_{\bH}(t,x,y)\,
\geq\, \int_{\sB(y,t)}\int_{\sB(x,t)}
p_{\bH}(t/3,x,u) p_{\bH}(t/3,u,v)p_{\bH}
(t/3,v,y)dudv \nonumber\\
\geq&  \left(\inf_{(u,v) \in \sB(x,t) \times
\sB(y,t)} p_{\bH}(t/3,u,v)\right)\int_{\sB(y,t)}
\int_{\sB(x,t)}
p_{{\bH}}(t/3,x,u)p_{\bH}(t/3,v,y)dudv.
\end{align*}
Thus \eqref{e:loww1} follows from Lemma \ref{l:loww_01}.

Observe that for $(u,v) \in \sB(x,t) \times \sB(y,t)$,
\begin{align}\label{e:deltage}
|\xi_x(t)-\xi_y(t)|=|x-y|, \qquad\delta_{\bH}(u) \wedge \delta_{\bH}(v) \ge
 a_1  \Phi^{-1}(t),
 \end{align}
and
 \begin{equation}\label{e:7.2}
 |x-y| -2M_1 \Phi^{-1}(t) \le |u-v| \le |x-y| +|u-\xi_x(t)|+|v-\xi_y(t)| \le |x-y| +2M_1 \Phi^{-1}(t).
 \end{equation}

When $|x-y|  > 4M_1 \Phi^{-1}(t)$, we have by \eqref{e:7.2} that for $(u,v) \in \sB(x,t) \times \sB(y,t)$,
$$
 |x-y| /2 \le  |u-v| \le 3|x-y|/2
$$
and so $2M_1 \Phi^{-1}(t) \le  |u-v|$.
Thus, for
$|x-y|  > 4M_1 \Phi^{-1}(t)$,
\begin{align}
\inf_{(u,v) \in \sB(x,t) \times \sB(y,t)} p_{\bH}(t/3,u,v) \ge
\inf_{(u,v) : 2M_1 \Phi^{-1}(t) \le  |u-v|  \le 3|x-y|/2\atop
\Phi(\delta_{\bH}(u)) \wedge \Phi(\delta_{\bH}(v)) >  a_1t} p_{\bH}(t/3,u,v).
 \label{e:flow11}
   \end{align}

 When $|x-y|  \le 4M_1 \Phi^{-1}(t)$, by \eqref{e:7.2} $|u-v|\le 6M_1 \Phi^{-1}(t)$ for $(u,v) \in \sB(x,t) \times \sB(y,t)$. Thus using  \eqref{e:deltage} and ${\bf PHI(\Phi)}$   (at most
$2+ 12[ M_1/a_1 ]$ times) and Lemma \ref{L:1.1} and Proposition \ref{step31}, we get
\begin{align}\label{e:fl11}
p_{\bH}(t/3,u,v) \ge c_1 p_{\bH}(t/6,u,u)  \ge c_{2} \Phi^{-1}(t) \quad \text{ for every } (u,v) \in \sB(x,t) \times \sB(y,t).
\end{align}
\eqref{e:loww21} now follows from \eqref{e:loww1}, \eqref{e:flow11} and \eqref{e:fl11}.
\qed

\section{Heat
kernel upper bound estimates in half spaces}\label{S:5}

In
this section, we consider a large class of symmetric L\'evy processes with concrete condition on the L\'evy densities.
Under these conditions, we can check that conditions
   {\rm ({\bf Comp})}, {\rm ({\bf HKC})} and
${\bf PHI(\Phi)}$ all hold. Thus we can apply
Proposition \ref{p:u_off1n} and Theorem \ref{t:prlower}
to
establish sharp two-sided estimates of the transition density of such L\'evy processes in half spaces.

Suppose that $\psi_1$ is an increasing function on
$[ 0, \infty )$ with $\psi_1(r)=1$ for $~0<r\leq 1$
 and there are
 constants $a_2\geq a_2>0$,
 $\gamma_2\ge \gamma_1> 0$ and $\beta \in [0,\infty]$ so that
\begin{equation}\label{eqn:exp}
a_1  e^{\gamma_1r^{\beta}} \leq \psi_1 (r)\leq a_2  e^{\gamma_2r^{\beta}}
\qquad \hbox{ for every }~ 1<r<\infty.
\end{equation}
 Suppose that
 $\phi_1$ is a strictly increasing function on $[0, \infty)$
with $\phi_1 (0)=0$, $\phi_1 (1)=1$ and
there exist constants
$0<a_3<a_4$ and $0<\beta_1\leq \beta_2<2$
so that
\begin{equation}\label{eqn:poly}
a_3 \Big(\frac Rr\Big)^{\beta_1} \,\leq\,
 \frac{\phi_1 (R)}{\phi_1 (r)}  \ \leq \ a_4
\Big(\frac Rr\Big)^{\beta_2}
\qquad \hbox{for every } 0<r<R<\infty.
\end{equation}
Since $0<\beta_1\leq \beta_2<2$, \eqref{eqn:poly} implies that
\begin{equation}\label{eqn:int1}
\int_0^r\frac {s}{\phi_1 (s)}\, ds \,\asymp \, \frac{ r^2}{\phi_1
(r)}, \qquad \int_r^\infty \frac{1}{s \phi_1 (s)}ds \ \asymp \, \frac{ 1}{\phi_1
(r)} \qquad   \hbox{for every } r>0.
\end{equation}

Throughout the remainder of this paper, we assume that
 {\rm({\bf UJS})} holds and that
  there are constants
$\gamma\geq 1$, $\kappa_1$, $\kappa_2$ and $a_0\geq 0$ such that
\begin{equation}\label{e:5.4}
\gamma^{-1} a_0 |\xi|^2 \leq \sum_{i, j=1}^d a_{i, j} \xi_i \xi_j \leq \gamma a_0 |\xi|^2 \quad \hbox{for every } \xi \in \R^d,
\end{equation}
and
\begin{equation}\label{e:5.5}
\gamma^{-1} \frac{1}{|x|^d \phi_1 (|x|)\psi_1 (\kappa_2|x|)} \\
 \le J(x) \le\gamma \frac{1}{|x|^d \phi_1 (|x|)\psi_1 (\kappa_1|x|)}\quad \hbox{for } x \in \R^d.
\end{equation}
 Note that ({\bf UJS}) holds  if $\kappa_1=\kappa_2$ in \eqref{e:5.5}.

Recall $\Phi $ is the function defined in \eqref{e:1.9}.
The next lemma gives explicit relation between $\Phi$ and $\phi_1$.

\begin{lemma}\label{L:5.1}
When $\beta = 0$,
\bee\label{e:5.3}
 \Phi ( r  ) \asymp \begin{cases}
 \phi_1 (r ) \1_{\{a_0=0\}} + r^2  \1_{\{a_0>0\}} \quad & \hbox{for } r \in [0, 1] , \\
 \phi_1 (r)  &\hbox{for } r\geq 1;
 \end{cases}
\eee
while for $\beta \in (0, \infty]$,
\begin{align}\label{e:phiprp}
\Phi(r) \asymp
\begin{cases}
\phi_1(r)\1_{\{a_0=0\}} + r^2 \1_{\{a_0>0\}} \quad  & \text{for  } r \in [0, 1], \\
r^2   & \text{for  } r \ge 1.
\end{cases}
\end{align}
\end{lemma}

\pf By Lemma \ref{l:j-upper} and \eqref{e:5.4},
$$
\frac1{\Phi(r)}\asymp \frac{a_0}{r^2}+\int_{\R^d} \left(1 \wedge \frac{|z|^2}{r^2}\right) J(z)dz
$$
Thus, by \eqref{e:5.5} and \eqref{eqn:exp}
\begin{align}\label{e:5.5b}
&c_0^{-1}\left(\frac{a_0}{r^2}+ r^{-2}\int_0^r
\frac{s}{\phi_1 (s)} e^{-\kappa_2 \gamma_2s^\beta }ds+\int_r^\infty \frac{1}{s \phi_1 (s)} e^{-\kappa_2 \gamma_2s^\beta }ds \right)\nn\\
& \le \frac1{\Phi(r)}\le  c_0 \left(\frac{a_0}{r^2}+ r^{-2}\int_0^r
\frac{s}{\phi_1 (s)} e^{-\kappa_1 \gamma_1s^\beta }ds+\int_r^\infty \frac{1}{s \phi_1 (s)} e^{-\kappa_1 \gamma_1s^\beta }ds \right).
\end{align}
When $\beta=0$, it follows from
  \eqref{eqn:int1} and \eqref{e:5.5b}  that
\begin{align}\label{e:beta=0}
\frac1{\Phi(r)}
\asymp \frac{a_0}{r^2}+  r^{-2}\int_0^r \frac{s}{ \phi_1 (s)}ds
+\int_r^\infty \frac{1}{s \phi_1 (s)}ds
 \asymp   \frac{a_0}{r^2}+ \frac{1}{\phi_1(r)} \qquad \hbox{for } r  >0.
\end{align}
Note that taking $R=1$ and $r=1$ in \eqref{eqn:poly}, we have
\bee \label{e:5.7}
 \phi_1 (r) \geq a_4^{-1} r^{\beta_2} \geq a_4^{-1}r^2
 \quad\hbox{for } r\in [0, 1] \quad \hbox{ and } \quad
 \phi_1 (R) \leq a_4 R^{\beta_2} \leq a_4 R^2 \quad\hbox{for } R\geq 1.
\eee
This together with \eqref{e:beta=0} establishes \eqref{e:5.3}.

When $r\geq 1$ and $\beta>0$,
$$
\int_r^\infty s^{-\beta_1-1} e^{-\kappa_1 \gamma_1 s^\beta }ds \le c_1 \int_r^\infty s^{-3} ds \le c_1 r^{-2}/2.
$$
Thus by \eqref{eqn:poly}, for  $\beta>0$ and $r\geq 1$,
\begin{align*}
&r^{-2}\int_0^r
\frac{s}{ \phi_1 (s)} e^{-\kappa_1 \gamma_1s^\beta }ds+\int_r^\infty \frac{1}{s \phi_1 (s)} e^{-\kappa_1 \gamma_1s^\beta }ds
\\
&\le r^{-2}\int_0^\infty
\frac{s}{ \phi_1 (s)} e^{-\kappa_1 \gamma_1s^\beta }ds+\int_r^\infty \frac{1}{s \phi_1 (s)} e^{-\kappa_1 \gamma_1s^\beta }ds  \le c_2 r^{-2},
\end{align*}
while
\begin{align*}
r^{-2}\int_0^r
\frac{s}{ \phi_1 (s)} e^{-\kappa_2 \gamma_2s^\beta }ds+\int_r^\infty \frac{1}{s \phi_1 (s)} e^{-\kappa_2 \gamma_2s^\beta }ds \ge r^{-2}\int_0^1
\frac{s}{ \phi_1 (s)} e^{-\kappa_2 \gamma_2s^\beta }ds \geq c_3 r^{-2} .
\end{align*}

By \eqref{eqn:int1},
for $r\le 1$ and  $\beta>0$,
\begin{eqnarray*}
&& r^{-2}\int_0^r
\frac{s}{ \phi_1 (s)}e^{-\kappa_1 \gamma_1s^\beta }ds
+\int_r^\infty \frac{1}{s \phi_1 (s)} e^{-\kappa_1 \gamma_1s^\beta }ds \\
&\le& r^{-2}\int_0^r
\frac{1}{s \phi_1 (s)} ds+\int_r^1 \frac{1}{s \phi_1 (s)} ds
+ \int_1^\infty   e^{-\kappa_1 \gamma_1s^\beta }ds   \\
& \le& \frac{c_4} {\phi_1(r)}+ c_4  \le \frac{c_5} {\phi_1(r)}
\end{eqnarray*}
and
\begin{align*}
r^{-2}\int_0^r \frac{s}{ \phi_1 (s)}e^{-\kappa_2 \gamma_2s^\beta }ds
+\int_r^\infty \frac{1}{s \phi_1 (s)} e^{-\kappa_2 \gamma_2s^\beta }ds
\ge e^{-\kappa_2 \gamma_2} r^{-2}\int_0^r \frac{s}{ \phi_1 (s)}ds
\ge \frac{c_6} {\phi_1(r)}.
\end{align*}
These combined with \eqref{e:5.5b} and \eqref{e:5.7}
immediately yield \eqref{e:phiprp}.    \qed

As an immediate consequence of Lemmas and \ref{l:Cnew} and \ref{L:5.1}, we have the following.

\begin{cor}\label{c:JUJS}
The conditions
\eqref{e:4.1},
  {\rm ({\bf ExpL})} and  {\rm ({\bf Comp})} hold.
\end{cor}

Since we have assumed  {\rm({\bf UJS})},
\eqref{e:5.4}
and \eqref{e:5.5},
our L\'evy process $X$ belongs to a subclass of the processes considered in
\cite{CK2, CK3, CKK3, CKK14}.
Therefore  $p(t, x, y)$ is  H\"older continuous  on $(0, \infty )\times \R^d\times \R^d$
and for every open set $D$, transition density $p_D(t, x, y)$ for the killed process $X^D$
  is  H\"older continuous on $(0, \infty )\times D\times D$ .
Define
\bee \label{e:5.8}
p^c (t, r)=t^{-d/2} \exp (- r^2/t).
\eee
Recall that $a_0$ is the ellipticity constant in \eqref{e:5.4}.
For each  $a, T>0$, we define a function $h_{ a,  T}(t, r)$ on $(t,r) \in (0, T] \times [0, \infty)$ as
\begin{equation}\label{eq:qd}
h_{a, T}(t, r):
=\begin{cases}
\displaystyle  a_0 p^c(t, ar)+
\big(\Phi^{-1}(t)^{-d} \wedge \big( tj(ar) \big)
 \quad &\hbox{if } \beta\in[0,1] \hbox{ or } r\in [0, 1],\\
 t \exp\left(-a \left( r \, \left(\log \frac{ T r}{t}\right)^{(\beta-1)/\beta}\wedge  r^\beta\right) \right)
 \qquad &\hbox{if } \beta\in(1, \infty) \text{ with } r \ge 1,\\
\displaystyle \left(t/(T r)\right)^{ar}
&\hbox{if } \beta=\infty \text{ with } r \ge 1; \\
\end{cases}
\end{equation}
and, for each  $a, T>0$,
define a function $k_{a,  T}(t, r)$ on $(t,r) \in [T, \infty) \times [0, \infty)$ as
\begin{equation}\label{eq:qd2}
k_{a, T}(t, r):
=\begin{cases}
\displaystyle   \Phi^{-1}(t)^{-d}  \wedge \left[ (a_0 p^c(t, ar)) +   tj(ar)  \right] \quad &\hbox{if } \beta =0,\\
t^{-d /2}
\exp\left( -a(r^\beta \wedge \frac{r^2 T}{t})  \right) &\hbox{if } \beta\in (0,1],\\
t^{-d /2}
\exp\left(-a r \left(      \big(1+ \log^+ \frac{rT}{t}\big)^{(\beta-1)/\beta} \wedge \frac{r  T}{t}\right) \right)  &\text{if } \beta\in (1, \infty),\\
t^{-d /2}
\exp\left(-a r \left(  \big(1+ \log^+ \frac{ rT }{t}\big)\wedge \frac{r^2T}{t}  \right)\right)
&\hbox{if }  \beta=\infty.
\end{cases}
\end{equation}
Note that
$r \to h_{a,  T}(t, r)$ and $r \to k_{a,  T}(t, r)$ are decreasing.

\bigskip

\begin{thm}\label{T:1.1}
The parabolic Harnack inequality ${\bf PHI(\Phi)}$ holds. Moreover,
for each positive constant $T$, there are positive constants
 $c_i$, $i=1, \dots 6$, which depend on the ellipticity constant $a_0$ of {\rm \eqref{e:5.4}},  such that
$$
c_2^{-1} h_{c_1,  T}(t, |x|)\, \le\, p(t,x)\,\le\,
c_2  h_{c_3,  T}(t, |x|) \quad \text{for every } (t, x) \in (0, T] \times \R^d,
$$
and
$$
c_4^{-1}k_{c_5, T}(t, |x|)\,\le\,
p(t,x)\,\le\,
c_4\,k_{c_6, T}(t, |x|) \quad \text{for every } (t, x) \in [T, \infty) \times \R^d.
$$
In particular, the condition {\rm ({\bf HKC})} holds.
\end{thm}

The above two-sided estimates on $p(t, x, y)$
follow  from \cite[Theorem 1.2]{CK2} and \cite[Theorems 1.2 and 1.4]{CKK3} when $a_0=0$, and from
\cite{CK3, CKK14} when $a_0>0$. Note that even though in \cite[Theorem 1.2]{CK2} and \cite[Theorems 1.2 and 1.4]{CKK3} two-sided estimates for $p(t, x, y)$ are stated separately for the cases $0<t\le 1$ and $t>1$, the constant 1 does not play any special role.
In fact, for example for $T<1$ one can easily check
 $$
c_3^{-1} h_{c_2,  1} (t, r)  \le h_{c_1,  T} (t, r) \le c_3 h_{c_2,  1} (t, r) \quad \text{ on  }t <T,
$$
and the  two-sided estimates for $
p(t,x)$ hold for the cases $0<t\le T$ and $t>T$, and can be stated in
the above way.

\begin{remark}\label{R:5.4} \rm
We remark here that in \cite[Theorems 1.2(2.b)]{CKK3},
the $| \log \frac{|x-y|}{t} |$ term should replaced by $1 + \log^+ \frac{|x-y|}{t}$. In the proof of \cite[Theorems 1.2(2.b)]{CKK3}, the case that $|x-y| \asymp t$ when $\beta\in (1, \infty)$ missed to be considered. Once taking into account of that missed case,  One can see from \cite{CKK3}  that \eqref{eq:qd2} is the correct form. See the statement and the proof of Proposition \ref{step6} below for the lower bound.
\end{remark}

We now present the main result of this section.

\begin{thm}\label{t:up3}
There exist
$c_1,  c_2>0$ such that for all $(t, x, y)\in (0, \infty) \times
\bH\times \bH$,
\begin{eqnarray*}
 p_{\bH}(t, x, y)&\le& c_1 \left(\sqrt{\frac{\Phi(\delta_{\bH}(x))}{t}} \wedge 1\right)\left(\sqrt{\frac{\Phi(\delta_{\bH}(y))}{t}} \wedge 1\right)
 \begin{cases}
  h_{c_2, 1}(t, |x-y|/6) & \text{ if } t \in (0, 1),\\
 k_{c_2, 1}(t, |x-y|/6) & \text{ if } t \in [1, \infty) .
  \end{cases}
\end{eqnarray*}
\end{thm}

\pf Since $r \to h_{a,  T}(t, r)$ and $r \to k_{a,  T}(t, r)$ are decreasing,  by Theorem \ref{T:1.1}
\begin{eqnarray*}
\sup_{w:  |w| \ge \frac{|x-y|}{6}} p(t, w)
 &\le& c_1
 \begin{cases}
\sup_{w:  \frac{|x-y|}{6}  \leq |w| }  h_{c_2, 1}(t, |w|) & \text{ if } t \in (0, 1),\\
\sup_{w:  \frac{|x-y|}{6}  \leq |w| }  k_{c_2, 1}(t, |w|) & \text{ if } t \in [1, \infty),
  \end{cases}
  \nn\\
 &\le& c_1
 \begin{cases} h_{c_2, 1}(t, |x-y|/6) & \text{ if } t \in (0, 1),\\
 k_{c_2, 1}(t, |x-y|/6) & \text{ if } t \in [1, \infty).
  \end{cases}
\end{eqnarray*}
This together with Proposition \ref{p:u_off1n} proves the theorem.
\qed

\section{ Interior lower bound estimates}\label{S:6}

In this section, we derive
following preliminary lower bound estimates on $p_\bH (t, x, y)$.
Recall that  we have assumed  {\rm({\bf UJS})},
\eqref{e:5.4} and \eqref{e:5.5}.

\begin{thm}\label{t:inlow}
Let  $a, T$ be positive constants.
There exist $c=c(a,
\beta_1, \beta_2, \beta, T)>0$ and $C_4=C_4(a,
\beta_1, \beta_2,  \beta, T)>0$ such that
$$
p_{\bH}(t, x, y) \ge c \begin{cases}
  h_{C_4, T}(t, |x-y|) & \text{ if } t \in (0, T),\\
 k_{C_4, T}(t, |x-y|)  & \text{ if } t \in [T, \infty),
  \end{cases}
$$
for every $(t, x, y)\in
(0, \infty)\times {\bH}\times {\bH}$ with
$ \delta_{\bH}(x) \wedge \delta_{\bH} (y) \ge a \Phi^{-1}(t)$.
\end{thm}

We will prove this theorem  through several propositions.
The following proposition
follows immediately from
Propositions~\ref{step1} and~\ref{step3}, Lemma \ref{L:5.1} and condition  \eqref{e:5.5}.

\begin{prop}\label{step31}
Let $D$ be an open subset of $\R^d$.
For every  $a>0$,
there exists a constant $c=c(a)>0$ so that
$$
p_D(t, x, y)\,\ge  \,c \, ((\Phi^{-1}(t))^{-d} \wedge  { t}{j(|x-y|)})
$$
for every $(t, x, y)\in
(0, \infty)\times D\times D$ with
$ \delta_D(x) \wedge \delta_D (y) \ge a \Phi^{-1}(t)$.
\end{prop}

Proposition \ref{step31} yields the interior lower bound for $p_D(t,x,y)$ and $p(t, x, y)$ for the case $\beta =0$ and $a_0=0$.
Proposition \ref{step31} also yield the interior lower bound for $p_D(t,x,y)$ and $p(t, x, y)$ for the case $\beta \in (0, 1]$, $t \le T$ and $a_0=0$.
As a direct  consequence of Proposition \ref{step3}, we have

\begin{cor}\label{C:6.2}
Suppose $\beta \in (0, \infty)$. For every
$a, T, C_*>0$,
there exist $c_1, c_2>0$
so  that
$$
p_{\bH}(t, x,y)\geq c_1
 \,  t\, e^{-c_2 |x-y|^\beta}
 $$
for every $(t, x, y)\in
[T, \infty)\times {\bH}\times {\bH}$ with
$ \delta_{\bH}(x) \wedge \delta_{\bH} (y) \ge a \Phi^{-1}(t)$ and $|x-y| \ge C_* \Phi^{-1}(t)$.
 In particular, when $0 <\beta \le 1$,
for every $a, T, C_*>0$, there exist $c_1,c_2>0$ such that
$$
p_{\bH}(t, x,y)\geq c_1
 \,  t\, e^{-c_2 |x-y|^\beta}\quad \hbox{when } |x-y|^{2-\beta}\geq t/C_*,  \,  \delta_{\bH}(x) \wedge \delta_{\bH} (y) \ge a \Phi^{-1}(t)
\hbox{ and } t\geq T.
$$
\end{cor}

The last assertion in Corollary \ref{C:6.2} holds because
$\Phi^{-1}(t) \asymp t^{1/2}$ for $t\geq T$ (by Lemma \ref{L:5.1}), and
for $t\geq T$ and $x, y$ with $|x-y|^{2-\beta} \geq t/C_*$,
one has $|x-y|^2 \geq ct$ where $c= (T/C_*)^{1/(2-\beta)} C_*^{-1}$.

A standard chaining argument give the following Gaussian lower bound. The proof is similar to the one of \cite[Theorem 5.4]{CKK3}.

\begin{prop}\label{main50}
Suppose
$\beta \in (0, \infty]$.
For every $C_*, a, T>0$, there exist constants $c_1,c_2>0$
 such that
$$
p_{\bH}(t, x,y)\geq c_1  t^{-d /2}\exp \left( -\frac{c_2 |x-y|^2}t \right)
$$
for every $(t, x, y)\in
[T, \infty)\times {\bH}\times {\bH}$ with
$ \delta_{\bH}(x) \wedge \delta_{\bH} (y) \ge a \Phi^{-1}(t)$ and $C_*|x-y| \le t/T$.
\end{prop}

\pf By considering $t/T$ instead of $t$, without loss of generality we assume
$T=1$.
Fix a constant $C_*>0$ and let $R:=|x-y|$.
When $t\geq 1\geq R$, by Proposition \ref{step31} and \eqref{e:phiprp},
$p_D(t, x, y) \geq c_1 \Phi^{-1} (t)^{-d} \geq c_2 t^{-d/2}$.
When  $t\geq  R^2\geq 1$, note that $R^2 \leq c_3 \Phi  ( r ) $ for some $c_3 >0$.
Thus in view of Lemma \ref{L:1.1},
by applying the parabolic Harnack inequality ${\bf PHI (\Phi)}$
at most $3(1+16a^{-2})c_3$ times, we have from  Proposition \ref{step1} that
$p_D(t, x, y) \geq c_4 \Phi^{-1} (t)^{-d} \geq c_5 t^{-d/2}$.
Hence we only need to consider the case
$ 1\vee (C_* R) \le t \le R^2$ (so $C_* \le 1$), which we now assume.
By \eqref{e:phiprp}, there exist a constant $c_0 \in (0,1)$ such that
$$
c_0^{-1} \sqrt s \ge \Phi^{-1}(s) \ge c_0 \sqrt s \qquad \hbox{for every }
s  \ge 2^{-1} (C_*)^2.
$$
Thus $\delta_{\bH}(x) \wedge \delta_{\bH} (y) \ge a  c_0\sqrt t$.

 Let $n$ be the smallest positive integer so that $t/n\ge (R/n)^2$.  Then
\bee \label{e:6.1}
1 \le R^2/t  \le n < 1+R^2/t  \le 2R^2/t \quad \hbox{ and } \quad
2(R/n)^2  \ge t/n\ge (R/n)^2.
\eee
Since $t\ge C_*R$, by \eqref{e:6.1}
\bee \label{e:6.1n}
\frac{t}{n} \ge \frac{t}{1+R^2/t}= \frac{t^2}{t+R^2} \ge 2^{-1} \left( \frac{t}{R}\right)^2 \ge  2^{-1} (C_*)^2.
\eee
 Let $x=x_0,x_1,\cdots,x_n=y$ be the points equally spaced on the line segment
 connecting $x$ to $y$ so that $|x_i-x_{i+1}| =R/n$ for
$i=0,\cdots, n-1$. Set  $B_i:=B(x_i,
2^{-1}a c_0R/n)$. Since
$t/n \geq (R/n)^2$ (by \eqref{e:6.1}) and $t/n\ge 2^{-1} (C_*)^2 $ (by \eqref{e:6.1n}),
we have for every   $(y_i,y_{i+1}) \in B_i \times B_{i+1}$,
$$\delta_{\bH} (y_i) \wedge \delta_{\bH} (y_{i+1}) \ge
a  c_0\sqrt t -2^{-1}ac_0R/n \ge 2^{-1} a  c_0\sqrt{t/n} \ge 2^{-1} a  c_0^2 \Phi^{-1}(t/n)$$
and
 $$ 4|y_i-y_{i+1}|
 \le 4 (1+2^{-1}ac_0) R/n   \le  4 (1+2^{-1}ac_0)\sqrt{t/n}  \le   4 (c_0^{-1}+2^{-1}a) \Phi^{-1}(t/n).$$
By Proposition \ref{step1} and applying ${\bf PHI (\Phi)}$ at most $N$ times,
where $N$ depends only on
$a$, $\Phi$ and $C_*$ (or by Proposition \ref{step31}),
we have
\begin{equation}
  p_{\bH} (t/n, y_i,y_{i+1})\geq c_6(t/n)^{-d/2}  \quad \text{for every } (y_i,y_{i+1}) \in B_i \times B_{i+1}.
\label{chalb22}\end{equation}
Using \eqref{chalb22} and then \eqref{e:6.1}, we have
\begin{eqnarray}
p_{\bH}(t, x,y)&\geq &
\int_{B_1}\dots\int_{B_{n-1}} p_{\bH} (t/n, x,y_1)\dots p_{\bH}
( t/n, y_{n-1},y)dy_1\dots dy_{n-1}\nn\\
&\geq & c_6(t/n)^{-d/2} \Pi_{i=1}^{n-1} \left(
c_7(t/n)^{-d/2}(R/n)^d \right)
\ge c_8(t/n)^{-d/2}(c_7 2^{-d/2})^{n-1}\nn\\
&\ge & c_8(t/n)^{-d/2}\exp (-c_9n)
\ge  c_{10}t^{-d/2}\exp \left( -\frac{c_{11}|x-y|^2}t \right). \label{e:newgb}
\end{eqnarray}
\qed

\begin{prop}\label{P:6.5}
Suppose $a_0, a>0$.
There are positive constants $c_1$ and $c_2$ so that
$$
p_{\bH}(t, x,y)\geq c_1 \Phi^{-1}(t)^{-d} \wedge
 \left( t^{-d /2}\exp \left( - \frac{c_2 |x-y|^2}t \right) + tj(|x-y|)\right)
$$
for every $(t, x, y)\in
[0, \infty)\times {\bH}\times {\bH}$ with
$ \delta_{\bH}(x) \wedge \delta_{\bH} (y)\geq a \Phi^{-1}(t) $
if either $\beta \in [0, 1]$ or $|x-y|\leq 1$.
\end{prop}

\pf  We first consider following five cases: (1)
$t\geq 1$ and $|x-y|\leq 1$ when $\beta \in [0, \infty]$, (2) $t\geq 1$ and $x, y \in \R^d$ when $\beta =0$,
(3) $|x-y|^{2-\beta} \ge t \ge 1$ when $\beta  \in (0, 1]$,
(4)  $t\leq 1$ and $|x-y|\geq 1$ when $\beta  \in (0, 1]$ (5) $|x-y|^2 \leq t \leq 1$ when $\beta \in [0, \infty]$.

Using the condition
\eqref{e:5.5}, we see  that for these five cases it holds that for every $c_1>0$ there is $c_2>0$ such that
\bee \label{e:6.3}
t^{-d /2}\exp \left( - \frac{c_1 |x-y|^2}t \right)\leq c_2 t j(|x-y|).
\eee
Hence by  Propositions \ref{step31} and
\ref{main50} and
\eqref{e:5.3}-\eqref{e:phiprp}, it suffices to consider the case
when $t\leq |x-y|^2 \leq 1$, which we will assume for the remainder of the proof.

By \eqref{e:5.3}-\eqref{e:phiprp}, there is a constant $c_1 \in (0, 1/2)$ so that
$c_1 r^2  \leq \Phi(r) \leq r^2/c_1 $. Set $R=|x-y|$.
Let $n$ to be the smallest integer so that $t/n \geq c_1^{-1} (R/n)^2$.
Observe that $ \frac{R^2}{c_1t} \leq n \leq \frac{2R^2}{c_1 t}$.
Let $x_0=x$, $x_1, \dots, x_n=y$ be the evenly spaced points on
be the line segment connecting $x$ to $y$ so that $|x_i- x_{i+1}|=R/n$.
Let $B_i=B_i(x_i, a R/4n)$. Then for every
we have for every   $(y_i,y_{i+1}) \in B_i \times B_{i+1}$,
$$
\delta_{\bH} (y_i) \wedge \delta_{\bH} (y_{i+1})
\ge a  \sqrt {c_1 t} - a R/2n \geq a  \sqrt {c_1 t} /2
\ge 2^{-1} a  c_1  \Phi^{-1}(t)$$
and
 $$
  4|y_i-y_{i+1}|   \le 4 (1+a) R/n   \le
    8\sqrt{c_1t/n}  \le  8 \Phi^{-1}(t/n).
 $$
By Proposition \ref{step1} and applying ${\bf PHI (\Phi)}$ at most finite many
times (or by Proposition \ref{step31}), we have
\begin{equation}\label{e:6.4}
  p_{\bH} (t/n, y_i,y_{i+1})\geq c_3(t/n)^{-d/2}  \quad \text{for every } (y_i,y_{i+1}) \in B_i \times B_{i+1}.
\end{equation}
Using \eqref{e:6.4} and then \eqref{e:6.1}, by the same argument in \eqref{e:newgb} we have
\begin{eqnarray*}
p_{\bH}(t, x,y)  \ge  c_4t^{-d/2}\exp \left( -\frac{c_5|x-y|^2}t \right).
\end{eqnarray*}
This together with Proposition \ref{step31} gives the desired lower bound
interior estimate for $t\leq |x-y|^2 \leq 1$.
This completes the proof of the proposition.
\qed

\medskip

Proposition \ref{step31}, Corollary \ref{C:6.2} and Propositions
\ref{main50}-\ref{P:6.5}  give the desired interior lower bound
stated in Theorem \ref{t:inlow} for $p_{\bH}(t, x, y)$
 when  $\beta \in [0,1]$.
We now consider the case $\beta=\infty$ and the case
$\beta\in(1, \infty)$ separately.

\begin{prop}\label{step5}
Suppose that $T, a>0$ and $\beta=\infty$.
Then, there exist constants $c_i=c_i(a, T)>0$, $i=1,2$, such that for any $x, y$ in ${\bH}$ with $\delta_{\bH}(x)\wedge\delta_{\bH}(y) \ge a \Phi^{-1}(t)$,  we have
\begin{equation}\label{ew218}
p_{\bH}(t, x,y)\geq
 c_1 \left( \frac t{T|x-y|} \right)^{c_2 |x-y|}
\quad \hbox{when } |x-y|\geq 1 \vee (t/T).
 \end{equation}
\end{prop}

\pf By considering $t/T$ instead of $t$, without loss of generality we assume $T=1$.
We let $R_1:=|x-y| \ge 1$.
We define $k$ as the integer satisfying $(4 \le ) 4  R_1\leq k <4  R_1+1<5 R_1$ and $r_t:= 2^{-1} a \Phi^{-1}(t)$.
 Let $x=x_0,x_1,\cdots,x_k=y$ be the points equally spaced on the line segment
 connecting $x$ to $y$ so that $|x_i-x_{i+1}| =R_1/k$ for
$i=0,\cdots, k-1$
and  $B_i:=B(x_i, r_t)$, with $i=0,1,2,\ldots,k$.
Then, $\delta_{\bH}(x_i)> 2r_t$ and $B_i=B(x_i, r_t) \subset B(x_i, 2 r_t)\subset {\bH}$, with $i=0,1,2,\ldots,k$.

Since $4  R_1\le k$, for each $y_i\in B_i$ we have
\begin{eqnarray}\label{e:stst1}
|y_i-y_{i+1}|\leq |y_i-x_i|+|x_i-x_{i+1}|+|x_{i+1}-y_{i+1}| \leq
\frac{1}{8}+\frac{R_1}{k}+ \frac{1}{8}<\frac{1}{2} .
\end{eqnarray}
Moreover, $\delta_{\bH}(y_i)\ge\delta_{\bH}(x_i)-|y_i-x_i |>r_t>r_{t/k}$ and $t/k \le R_1/k \le 1/4$

Thus, by Proposition~\ref{step31} and~\eqref{e:stst1}, there are constants $c_i=c_i(  a)>0$, $i=1,2$, such that for $(y_i, y_{i+1})\in B_i\times B_{i+1}$ we have
\begin{align}
p_{\bH}(t/k, y_i, y_{i+1})
&\ge c_1\left(\phi_1^{-1}(t/k)^{-d} \wedge \frac{t/k}{|y_i-y_{i+1}|^{d} \phi_1(|y_i-y_{i+1}|)}\right)\nn \\
&\ge c_2\left(\phi_1^{-1}(t/k)^{-d} \wedge t/k \right)=c_2\left(1 \wedge (t/k) \right) =  c_2 \,t/k.\label{e:stst2}
\end{align}
Observe that $4R_1\le k < 2(k-1) < 8R_1$, $\phi_1^{-1}(t/k) \ge a_3^{1/\beta_1} (t/k)^{1/\beta_1}$  and $r_t\ge  r_{t/k}$.
Thus, from \eqref{e:stst2} we obtain
\begin{align*}
p_{\bH}(t,x,y)& \ge \int_{B_1}\ldots\int_{B_{k-1}}p_{\bH}(t/k,x,y_1)\ldots p_{\bH}(t/k, y_{k-1},y) dy_{k-1}\ldots dy_1\\
&\ge (c_2t/k)^k\Pi^{k-1}_{i=1} |B_i|  \ge (c_2t/k)^k c_3^{k-1}(t/k)^{d(k-1)/\beta_1}\\
&\ge  c_4 (c_5t/k)^{c_6k}  \ge  c_7 (c_8t/R_1)^{c_9R_1}  \ge c_{10} (t/R_1)^{c_{11}R_1}.
\end{align*}
 \qed

\begin{prop}\label{step6}
Suppose that $T>0$,
$a >0$  and $\beta\in(1, \infty)$.
Then, there exist constants  $c_i=c_i(a, \beta,T)>0$, $i=1,2$ such that for any $x, y$ in ${\bH}$ with $\delta_{\bH}(x)\wedge\delta_{\bH}(y) \ge a \Phi^{-1}(t)$  we have
\begin{eqnarray*}
p_{\bH}(t,x,y)\ge c_1 t \exp\left( -c_2 \left(|x-y|\left(\log\frac{T|x-y|}{t}\right)^{\frac{\beta -1}{\beta}}\wedge (|x-y|)^{\beta}\right)\right) \quad \text{if } t \le T, |x-y|>1,
\end{eqnarray*}
and
\begin{eqnarray*}
p_{\bH}(t,x,y)\ge c_1 t^{-d/2} \exp\left( -c_2 \left(|x-y|\left(1+\log^+\frac{T|x-y|}{t}\right)^{\frac{\beta -1}{\beta}}\right)\right) \quad \text{if } t > T, |x-y| >t/T.
\end{eqnarray*}
\end{prop}

\pf Without loss of generality we assume $T=1$.
We fix  $a >0$, and we let $R_1:=|x-y|$.

\noindent
(i) If $1 \le R_1 \le 3$ and $t \le 1$,  the proposition holds by virtue of Proposition~\ref{step31}.

\noindent
(ii) If $R_1(\log(R_1/t))^{(\beta-1)/\beta}\ge (R_1)^{\beta}$ (when $t \le 1$), the proposition holds also by virtue of Proposition~\ref{step31}.

\noindent
(iii) If $t > 1$  and $3t \ge R_1\ge t$, the proposition holds  by virtue of Proposition~\ref{main50}.

\noindent
(iv) We now assume $(t, R_1) \in ((0,1] \times (3, \infty)) \cup  ((1, \infty) \times (3t, \infty))$ and
 $R_1(\log(R_1/t))^{(\beta-1)/\beta}< (R_1)^{\beta}$, which is equivalent to $R_1 \exp\{-\left(R_1\right)^{\beta}\}< t$.
Note that $R_1/t>3$.

Let $k\ge 2$ be a positive integer such that
\begin{align}\label{e:k}
1<R_1 \left(\log\frac{R_1}{t}\right)^{-1/\beta}\le k<R_1 \left(\log\frac{R_1}{t}\right)^{-1/\beta}+1< 2R_1 \left(\log\frac{R_1}{t}\right)^{-1/\beta}.
\end{align}
We define
$r_t:=(2^{-1}a\Phi^{-1}(t/R_1))\wedge ((6)^{-1}(\log ({R_1}/{t}))^{1/\beta})$.
Then, by \eqref{e:k}  we have
\begin{align}
\left(2^{-1}a\Phi^{-1}(t/R_1)\right)\wedge \frac{R_1}{6 k} \le r_t\le\frac{1}{6}\left(\log\frac{R_1}{t} \right)^{1/\beta}< \frac{R_1}{3 k}.\label{e:r_t.2}
\end{align}
 Let $x=x_0,x_1,\cdots,x_k=y$ be the points equally spaced on the line segment
 connecting $x$ to $y$ so that $|x_i-x_{i+1}| =R_1/k$ for
$i=0,\cdots, k-1$
and  $B_i:=B(x_i, r_t)$, with $i=0,1,2,\ldots,k$.
Then, $\delta_{\bH}(y_i)\ge 2^{-1}  a \Phi^{-1}(t)>2^{-1}  a\Phi^{-1}(t/k)$ for every  $y_i \in B_i$.
Note that from \eqref{e:r_t.2} we obtain
\begin{align}\label{e:y_i}
\frac{1}{3 } \frac{R_1}{k} \le |x_i-x_{i+1}|-2r_t \le |y_i-y_{i+1}|\le|x_i-x_{i+1}|+2r_t \le\frac{5}{3 }\frac{R_1}{k}
\end{align} for every  $(y_i, y_{i+1}) \in B_i \times  B_{i+1}$.
We also observe that, by \eqref{e:k}
$$
\frac{t}{k} \le \frac{t}{R_1} (\log (R_1/t))^{1/\beta}
\le \sup_{s \ge 3} s^{-1} (\log s)^{1/\beta}  <\infty$$
and
$$
\frac{R_1}{2k} \ge \frac14 (\log (R_1/t))^{1/\beta} \ge
\frac12  (\log 3)^{1/\beta} >0.$$
Thus, using Proposition~\ref{step31} along with \eqref{e:k} and  \eqref{e:y_i} we obtain
\begin{align}\label{p(t/k)}
&p_{\bH}(t/k, y_i, y_{i+1}) \ge c_{1}\frac{t}{k} j(|y_i-y_{i+1}|)\ge \, c_{2}\frac{t}{k}\left(R_1/k\right)^{-d-\beta_2}e^{-c_{3}(R_1/k)^{\beta}}\nn\\
& \ge \,c_{4} \frac{t}{R_1}\left(\frac{k}{2R_1}\right)^{d+\beta_2-1}e^{-c_{3}(R_1/k)^{\beta}}\, \ge \, c_{4} \frac{t}{R_1}\left(\log \frac{R_1}{t}\right)^{-\frac{d+\beta_2-1}{\beta}}
\left(\frac{t}{R_1}\right)^{c_{3}} \,\ge \, c_{4}\left(\frac{t}{R_1}\right)^{c_{5}}.
\end{align}
Since the definition of $r_t$  yields
$$r_t \ge c_{6}\left((t/R_1)^{(\beta_2 \wedge \beta)^{-1}}  \wedge (\log ({R_1}/{t}))^{1/\beta} \right) \ge
  c_{7}(t/R_1)^{(\beta_2 \wedge \beta)^{-1}} .$$
 by using \eqref{e:k}, \eqref{p(t/k)} and the semigroup property we conclude that
\begin{align*}
p_{\bH}(t, x, y)&\ge \int_{B_1}\cdots\int_{B_{k-1}} p_{\bH}(t/k, x, y_1)\cdots p_{\bH}(t/k, y_{k-1}, y ) dy_1\cdots dy_{k-1}\\
&\ge
 c_{4}^{k}c_{7}^{k-1}\left(\frac{t}{R_1}\right)^{c_{5}k+(\beta_2 \wedge \beta)^{-1}(k-1)} \\
&\ge c_{8}\exp \left( -c_{9}k\log({R_1/t})\right) \\
&\ge c_{8}\exp\left( -c_{9}\left(R_1 \log\left( R_1/t\right)^{-1/\beta}+1\right)\log (R_1/t) \right)\\
&\ge c_{8} \exp\left( -2c_{9}\left(R_1 \log\left(R_1/t\right)^{\frac{\beta -1}{\beta}}\right)\right)\\
&\ge c_{8}
\begin{cases} \displaystyle
t \exp\left( -2c_{9}\left(R_1 \log\left( R_1/ t \right)^{\frac{\beta -1}{\beta}}\right)\right)  & \text{ if } (t, R_1) \in (0,1] \times (3, \infty)\\ \displaystyle
t^{-d/2} \exp\left( -2c_{9}\left(R_1\Big( 1+ \log^+\left( R_1/t \right)^{\frac{\beta -1}{\beta}}\Big)\right)\right) & \text{ if }(t, R_1) \in  (1, \infty) \times (3t, \infty).
  \end{cases}
\end{align*}
\qed

Propositions \ref{step5}-\ref{step6} together with
Proposition \ref{step31} and Propositions \ref{main50}-\ref{P:6.5}
yield the interior lower bound estimates of  Theorem \ref{t:inlow}
for $\beta \in (1, \infty]$.

\medskip

\begin{remark}\label{R:6.7}  \rm
Assume that $D$ is an connected open set with the following property: there exist $\lambda_1 \in [1, \infty)$ and $\lambda_2 \in (0, 1]$ such that for every $r \le 1$ and  $x,y$ in  $D$ with $\delta_D(x)\wedge \delta_D(y)\ge r$ there exists in $D$ a length parameterized rectifiable curve $l$ connecting $x$ to $y$ with the length $|l|$ of $l$
 less than or equal to $\lambda_1|x-y|$ and $\delta_D(l(u))\geq\lambda_2 r$ for $u\in[0,|l|].$

Under this assumption, we can also prove
Theorem \ref{t:inlow} on such $D$ with minor modifications. We omit the details here; see \cite[Section 3]{KK} for the case $t<T$ and $\phi(r)=r^\alpha$.
\end{remark}

\section{Two-sided heat kernel estimates}\label{S:7}

In this section we prove the two-sided estimates of $p_{\bH}(t, x, y)$ under
conditions \eqref{e:5.4} and  \eqref{e:5.5}.

\begin{thm}\label{t:final}
Suppose {\rm({\bf UJS})},
\eqref{e:5.4} and \eqref{e:5.5} hold.
There exist
$c_1, c_2, c_3>0$ such that for all $(t, x, y)\in (0, \infty) \times
\bH\times \bH$,
\begin{eqnarray*}
 p_{\bH}(t, x, y)\le c_1 \left(\sqrt{\frac{\Phi(\delta_{\bH}(x))}{t}} \wedge 1\right)\left(\sqrt{\frac{\Phi(\delta_{\bH}(y))}{t}} \wedge 1\right) \times
 \begin{cases}
  h_{c_2, 1}(t, |x-y|/6) & \text{ if } t \in (0, 1),\\
 k_{c_2, 1}(t, |x-y|/6) & \text{ if } t \in [1, \infty),
  \end{cases}
\end{eqnarray*}
and
\begin{eqnarray*}
 p_{\bH}(t, x, y)\ge c_1^{-1} \left(\sqrt{\frac{\Phi(\delta_{\bH}(x))}{t}} \wedge 1\right)\left(\sqrt{\frac{\Phi(\delta_{\bH}(y))}{t}} \wedge 1\right) \times
 \begin{cases}
  h_{c_3, 1}(t, 3|x-y|/2) & \text{ if } t \in (0, 1),\\
 k_{c_3, 1}(t, 3|x-y|/2)  & \text{ if } t \in [1, \infty) .
  \end{cases}
\end{eqnarray*}
\end{thm}

\pf By Theorem \ref{t:up3}, we only need to show the lower bound of $p_{\bH}(t, x, y)$.
For this, we will apply Theorem \ref{t:prlower}.
Let  $C_4$ be the constant $C_4$ in Theorem \ref{t:inlow} with $T=1/3$.

Since
$r \to h_{a,  T}(t, r)$ and $r \to k_{a,  T}(t, r)$ are decreasing,
we have by  Theorem \ref{t:inlow} that for
$|x-y|  > 4M_1 \Phi^{-1}(t)$,
\begin{align}
 \inf_{(u,v) : 2M_1 \Phi^{-1}(t) \le  |u-v|  \le 3|x-y|/2
\atop
\Phi(\delta_{\bH}(u)) \wedge \Phi(\delta_{\bH}(v)) >  a_1t
}p_{\bH}(t/3,u,v)  &\ge c_1 \begin{cases}
  h_{C_4, 1/3}(t/3, 3|x-y|/2) & \text{ if } t \in (0, 1),\\
 k_{C_4, 1/3}(t/3, 3|x-y|/2) & \text{ if } t \in [1, \infty),
  \end{cases}\nn\\
& \ge c_2 \begin{cases}
  h_{c_3, 1}(t, 3|x-y|/2) & \text{ if } t \in (0, 1),\\
 k_{c_3, 1}(t, 3|x-y|/2) & \text{ if } t \in [1, \infty).
  \end{cases} \label{e:flow1}
   \end{align}

When $6M_1 \Phi^{-1}(t) \ge r $ and $t \ge 1$,  by
\eqref{e:phiprp}, we have  $c_8 M_1 t^{1/2} \ge r $. Thus on $6M_1 \Phi^{-1}(t) \ge r $ and $t \ge 1$
\begin{align}\label{e:fl12}
k_{C_4, 1/3}(t/3, r) &\ge c_4 \begin{cases}
 \Phi^{-1}(t) ^{-d} &\text{ if } \beta =0,\\
t^{-d /2}
\exp\left( -C_4(r^\beta \wedge 3\frac{r^2}{t})  \right) &\text{ if } \beta\in (0,1],\\
t^{-d /2}
\exp\left( -C_4 \left( \Big(r \, \big(1+ \log^+ \frac{  r}{t}\big)^{\frac{\beta-1}{\beta}}\Big)\wedge 3\frac{r^2}{t}\right) \right)
  &\text{ if } \beta\in (1, \infty),\nn \\
t^{-d /2}
\exp\left( -C_4 \left( \Big(r \, \big(1+ \log^+ \frac{  r}{t}\big)\Big)\wedge 3\frac{r^2}{t}\right) \right)
&\hbox{ if }  \beta=\infty,
\end{cases}\\
&\ge c_{5}
 \Phi^{-1}(t) ^{-d}. \end{align}
 So by \eqref{e:fl12} and Theorem \ref{t:inlow}, for $|x-y|  \le 4M_1 \Phi^{-1}(t)$,
\begin{align}\label{e:ge4M_0}
\Phi^{-1}(t)
  \ge c_{6}\begin{cases}
  h_{c_7, 1}(t, 3|x-y|/2) & \text{ if } t \in (0, 1),\\
 k_{c_7, 1}(t, 3|x-y|/2) & \text{ if } t \in [1, \infty).
  \end{cases}
\end{align}
Combining \eqref{e:loww21}, \eqref{e:flow1} and \eqref{e:ge4M_0}, we conclude that \begin{eqnarray*}
p_{\bH}(t, x, y)\ge c_{8} \left(\sqrt{\frac{\Phi(\delta_{\bH}(x))}{t}} \wedge 1\right)
\left(\sqrt{\frac{\Phi(\delta_{\bH}(y))}{t}} \wedge 1\right)\begin{cases}
  h_{c_9, 1}(t, 3|x-y|/2) & \text{ if } t \in (0, 1),\\
 k_{c_9, 1}(t, 3|x-y|/2) & \text{ if } t \in [1, \infty).
  \end{cases}
\end{eqnarray*}
\qed

\begin{remark}\label{R:7.2}  \rm
(i) In view of Theorem \ref{T:1.1}, we can restate Theorem \ref{t:final}
as follows. There are positive constants $c_i$, $1\leq i\leq 5$, so that
\begin{eqnarray}
&& \hskip -0.6truein \frac1{c_1} \left(\frac1{\sqrt{t\, \Psi (1/\delta_D(x))}}
\wedge 1 \right)  \left(\frac1{\sqrt{t\, \Psi (1/\delta_D(y))}} \wedge 1 \right)
 p(c_2t, c_3(y-x)) \,\le\, p_D(t, x,  y) \nonumber \\
&& \le\, c_1 \, \left(\frac1{\sqrt{t\, \Psi (1/\delta_D(x))}}
\wedge 1 \right)  \left(\frac1{\sqrt{t\, \Psi (1/\delta_D(y))}} \wedge 1\right)
 \, p(c_4t, c_5(y-x)) \label{e:7.4}
\end{eqnarray}
where  $p(t, x)$ is the transition density of $X$.
This essentially confirms the conjecture \eqref{e:1.10} for this class of symmetric L\'evy processes and for
$D=\bH$.

(ii) Recently sharp two-sided Dirichlet heat kernel estimates have been
established in \cite{CKS9, CKS6}
for a large class of symmetric L\'evy processes  in $C^{1,1}$ open sets
for $t\leq 1$. The L\'evy process considered in
\cite{CKS9, CKS6}
satisfy the conditions \eqref{e:5.4},  \eqref{e:5.5} and {\rm({\bf UJS})} of this paper.
Now assume $X$ is a symmetric L\'evy process considered \cite{CKS9, CKS6}.
Then using the
``push inward" method of
\cite{CT} (see \cite{CKS4} for its use in  relativistic stable processes case)
and the short time
heat kernel estimates in \cite{CKS9, CKS6},
 we can obtain
global sharp two-sided Dirichlet heat kernel estimates on
half-space-like $C^{1,1}$ open sets from the
the Dirichlet heat kernel estimates established in this paper
on half-spaces. We leave the details to the interested reader.
\end{remark}

{\bf Acknowledgement.} We thank the referees for helpful comments.

\vskip 0.3truein

{\bf Zhen-Qing Chen}

Department of Mathematics, University of Washington, Seattle,
WA 98195, USA

E-mail: \texttt{zqchen@uw.edu}

\bigskip

{\bf Panki Kim}

Department of Mathematical Sciences and Research Institute of Mathematics,

Seoul National University, Building 27, 1 Gwanak-ro, Gwanak-gu Seoul 08826, Republic of Korea

E-mail: \texttt{pkim@snu.ac.kr}


\vspace{.1in}

\small
\begin{thebibliography}{99}

\bibitem{BBCK}   {M.~T. Barlow,  R.~F. Bass,  Z.-Q. Chen and M. Kassmann,}
Non-local Dirichlet forms and symmetric jump processes.
\textit{Trans. Amer. Math. Soc.} \textbf{361} (2009), 1963--1999.

\bibitem{BG} K. Bogdan and T. Grzywny, Heat kernel of fractional Laplacian in cones. {\it Colloq. Math.} {\bf 118} (2010), 365--377.

\bibitem{BGR} K. Bogdan, T. Grzywny and M. Ryznar:
Heat kernel estimates for the fractional Laplacian with Dirichlet
Conditions. {\it  Ann. Probab.}
{\bf 38} (2010), 1901--1923.

\bibitem{BGR1} K. Bogdan, T. Grzywny and M. Ryznar,
Density and tails of unimodal convolution semigroups.
{\em J. Funct. Anal.} {\bf 266} (2014), 3543--3571.

\bibitem{BGR2} K. Bogdan, T. Grzywny and M. Ryznar,
Barriers, exit time and survival probability for unimodal L\'evy processes.
{\it Probab. Theory Relat. Fields} {\bf 162} (2015), 155--198.

\bibitem{BGR3} K. Bogdan, T. Grzywny and M. Ryznar,
Dirichlet heat kernel for unimodal L\'evy processes.
\emph{Stoch. Proc. Appl},  {\bf  124(11)} (2014), 3612--3650.

\bibitem{CKK2}  Z.-Q. Chen, P. Kim and  T. Kumagai,
On heat kernel estimates and parabolic Harnack inequality for jump processes on metric measure spaces.
\textit{ Acta Math. Sin. (Engl. Ser.)} \textbf{25} (2009), 1067--1086.

\bibitem{CKK3}  Z.-Q. Chen, P. Kim and  T. Kumagai,
Global heat kernel estimates for symmetric jump processes.
{\it  Trans. Amer. Math. Soc.} {\bf 363} (2011), 5021--5055.

\bibitem{CKK14}  Z.-Q. Chen, P. Kim and  T. Kumagai,
Heat kernel estimates for reflected diffusions with jumps on metric
measure spaces. In preparation.

\bibitem{CKS}  Z.-Q. Chen, P. Kim, and R. Song,
Heat kernel estimates for Dirichlet fractional Laplacian.
{\em J. European Math. Soc.} {\bf 12} (2010), 1307--1329.

\bibitem{CKS5} Z.-Q. Chen, P. Kim and R. Song,
Heat kernel estimate for $\Delta+\Delta^{\alpha/2}$ in $C^{1,1}$
open sets.
{\it J. London Math. Soc.}
{\bf 84}  (2011), 58--80.

\bibitem{CKSa1}  Z.-Q. Chen, P. Kim and R. Song,
Sharp heat kernel estimates for relativistic stable processes in open sets.
{\it Ann. Probab. \bf 40} (2012), 213-244.


\bibitem{CKS4}  Z.-Q. Chen, P. Kim, and R. Song,
Global heat kernel estimates for relativistic stable processes in
half-space-like open sets.
{\it Potential Anal.},
{\bf 36} (2012),  235--261.



\bibitem{CKS7}  Z.-Q. Chen, P. Kim, and R. Song:
Global heat kernel estimates for relativistic stable processes in exterior open sets.
{\it J. Funct. Anal.} {\bf 263} (2012), 448--475.


\bibitem{CKS9} Z.-Q. Chen, P. Kim and R. Song,
Dirichlet heat kernel estimates for rotationally symmetric L\'evy processes.
{\it Proc. London Math. Soc.} {\bf 109} (2014), 90--120.

\bibitem{CKS6} Z.-Q. Chen, P. Kim and R. Song,
Dirichlet heat kernel estimates for subordinate Brownian motions with Gaussian components.
 {\it J. Reine Angew. Math.}
 {\bf 711} (2016), 111--138.


 \bibitem{CK1} Z.-Q. Chen and  T. Kumagai,
 Heat kernel estimates for stable-like processes on $d$-sets. {\it Stochastic Process. Appl.}, {\bf 108} (2003), 27--62.



\bibitem{CK2} Z.-Q. Chen and  T. Kumagai,
Heat kernel estimates for jump processes of mixed types on metric measure spaces.
{\it Probab. Theory Relat. Fields}, {\bf 140} (2008), 277--317.

\bibitem{CK3}  Z.-Q. Chen and  T. Kumagai,
 A priori H\"older estimate, parabolic Harnack principle and heat kernel
estimates for diffusions with jumps.
{\it    Rev. Mat. Iberoam.} {\bf 26}
   (2010), 551--589.

\bibitem{CT} Z.-Q. Chen and  J. Tokle,
Global heat kernel estimates for fractional Laplacians in unbounded open sets.
{\it Probab. Theory Relat. Fields},
{\bf 149} (2011), 373--395.


\bibitem{D2}
 E. B. Davies,
The equivalence of certain heat kernel and Green function bounds.
{\em J. Funct. Anal.} {\bf 71} (1987), 88--103.

\bibitem{Dyn} E. B. Dynkin, {\it Markov processes, Vol. I}.
Academic Press, New York, 1965

\bibitem{G}
T. Grzywny,
 On Harnack inequality and H\"older regularity for isotropic unimodal L\'evy processes.
  {\it Potential Anal. \bf 41} (2014), 1-29.

\bibitem{KK}
P.~Kim and K.-Y. Kim, Two-sided estimates for the transition densities of symmetric Markov processes
dominated by stable-like processes in $C^{1,\eta}$ open sets.
\emph{Stoch. Proc. Appl},
{\bf124} (2014),  3055--3083.

\bibitem{KSV7} P. Kim, R. Song and Z. Vondra\v{c}ek, Uniform boundary Harnack principle for rotationally
 symmetric L\'evy processes in general open sets.
  {\it Sci. China Math.} {\bf 55}, (2012), 2193--2416.


 \bibitem{KS}
V. Knopova and R. Schilling,
A note on the existence of transition probability densities of L\'evy processes. {\it Forum Math.} {\bf 25} (2013), 125--149.

\bibitem{KMR} M. Kwa\'snicki, J. Ma{\l}ecki, M. Ryznar,  Suprema of L\'evy processes.
 {\em Ann.~Probab.} {\bf 41} (2013), 2047--2065.

 \bibitem{P}
W. E. Pruitt,
The growth of random walks and L\'evy processes.
{\it Ann. Probab.} {\bf 9} (1981), 948--956.

\bibitem{Sk}
A. V. Skorohod.
{\em Random Processes with Independent Increments}.
Kluwer, Dordrecht, 1991.


\bibitem{Zq3} Q. S. Zhang, The boundary behavior of heat kernels of Dirichlet Laplacians.
 {\em J. Differential Equations}, {\bf 182} (2002),
416--430.

\end{thebibliography}
\end{document}